\documentclass[10pt]{article}
\usepackage{graphicx}
\usepackage{amsmath}
\usepackage{amsfonts}
\usepackage{amssymb}

\begin{document}

\title{\textbf{\textquotedblleft Worthwhile-to-move\textquotedblright\
behaviors as temporary satisficing without too many sacrificing processes}}
\author{Hedy Attouch (ACSIOM-I3M)\thanks{%
Hedy Attouch, I3M UMR CNRS 5149, Universit\'e Montpellier II, Place Eug\`ene
Bataillon, 34095, Montpellier, France, email: attouch@math.univ-montp2.fr,
supported by French ANR program OSSDAA, 
ANR-08-BLAN-0294-03} \and Antoine Soubeyran (GREQAM) \thanks{%
Antoine Soubeyran, GREQAM UMR CNRS 6579, Universit\'e de la M\'editerran\'ee, 
Chateau Lafarge, Routes des Milles, 13290, Les
Milles, France, email: antoine.soubeyran@univmed.fr}}
\date{ }
\maketitle{}

\begin{abstract}
In a world with frictions, the \textquotedblleft
worthwhile-to-move\textquotedblright\ incremental principle is a mechanism
where, at each step, the agent, before moving and after exploration around
the current state, compares intermediate advantages and costs to change to
advantages and costs to stay. Theses advantages and costs to change are not
the usual benefit and cost functions. They are behavioral, including goal
setting, psychological, cognitive (learning) and inertia aspects. The agent
is supposed to have a long-term goal and limited needs. At each step, the
agent chooses between repeating the action or  changing it. Acceptable moves are such
that \textquotedblleft advantages to move than to stay,\textquotedblright\
are higher than some fraction of \textquotedblleft costs to move than to
stay,\textquotedblright \ with, as a result, a limitation of the intermediate
sacrifices to reach the goal. The transition process is made of a
punctuated succession of static exploration-exploitation phases, and dynamic
moving phases. This \textquotedblleft dynamic and reference
dependent\textquotedblright\ incremental cost-benefit behavior improves and
leads to local actions in an endogenous way. It converges if the agent has
high enough local costs to move, because of an entropy property.

When the agent is more goal-oriented and wants to \textquotedblleft improve
enough\textquotedblright \ at each step, the process shrinks (the diameter
of the state-dependent \textquotedblleft
worthwhile-to-move\textquotedblright \ set decreases to zero), and reduces
to one state if the goal function is upper semicontinuous. The process ends
in a permanent routine, a rest point (which may be lower than any local maximum) where
the agent prefers to stay than to change, in spite of some possible residual
frustration to have missed his goal. The convergence in a finite
number of steps occurs if the agent chooses a finite total exploitation time.

The model describes a full range of behaviors. Psychology and cognition play
an important role in the balance between motivation and fear to change.
Links with optimization theory, variational inequalities (Ekeland's
variational principle) are given. The \textquotedblleft as if
hypothesis\textquotedblright\ (as if agents would optimize) is revisited,
leading to the introduction of a new class of optimization algorithms
with inertia, called the \textquotedblleft local search and proximal
algorithms\textquotedblright.
\end{abstract}

\smallskip

\section{Introduction: From Incrementalism to Instrumentalism to Optimization%
}

In a world with frictions, the \textquotedblleft worthwhile-to-move
\textquotedblright\ incremental principle allows us to handle the extreme
but very realistic case of a \textquotedblleft muddling through" decision-making process 
where the agent at each step tries to do a little better
than before (Lindblom, 1959). By adding an instrumental goal-setting
principle,  the \textquotedblleft improving enough\textquotedblright\
principle, this model handles a great variety of intermediate goal-setting
behaviors, allowing us to use inertia and explaining the formation of
routines. Satisficing (Simon, 1956) is the most noteworthy example. The
traditional case of global optimization is at the other extreme of this
range of behaviors, in a world deprived of friction.

The term \textquotedblleft agent'' is taken in a broad sense, it may
represent a single  agent, a problem solver, or a group 
 whose structure remains identical during the process. We do not
examine here the case of interacting agents with strategic features.

Behavioral imperfections are the rule, not the exception. This is true for
human behaviors, for decision making, and for problems solving. The world is
full of frictions. Human beings have limited physical resources, bounded
cognitive abilities, psychological biaises which distort their evaluation
and their motivations to search. Human preferences and beliefs may be
inconsistent; motivations to act can be unconscious, vague, too low, or
impulsive; goals are frequently ill-defined; knowledge may be inadequate.
Agents cannot take advantage of any piece of information they receive. They
ignore a part of them. These imperfections generate costs to try to remedy.
Group members have non congruent goals and partial conflicts of interest,
they lose time to bargain and resolve conflicts. One can distinguish three
kinds of imperfection:

i) Cognitive imperfections: lack of knowledge is one limit for decision
making;

ii) Psychological imperfections: they are the main source of problems in
driving action properly. They concern psychological biais of evaluation,
editing, goal setting, motivation building, dealing with frustration
feelings, unclear goals, impulsive behaviors, emotion biais. We can mention
 reference-dependent biais, anchoring effects,
framing effects, bracketing, narrow or broad mental accounting, salient
effects, hyperbolic discounting, impatience, regret, negative feelings for
ambiguity, closure effects, inhibition and excitation, the role of stress
and emotions;

iii) Physical and physiological imperfections: inertia and frictions
(Rumelt, 1990) are the main problems for action.

Lack of knowledge was well documented, but goal setting and inertia,
frictions, and costs to change not as much. We first model the
imperfections related to inertia and frictions.

Because imperfections are everywhere, most of our behaviors are incremental.
They work step by step, using small improving steps and local improving
actions, local exploration devices, trials and errors. Incremental behaviors
characterize low goal-oriented \textquotedblleft
decision-making\textquotedblright\, where, at each step, the agent compares
local advantages  and  costs to move to advantages and costs to
stay. These incremental behaviors improve in time. Practical examples of behavioral
\textquotedblleft local cost-benefit processes" are \textquotedblleft pros
and cons lists\textquotedblright\ and \textquotedblleft plus-minus
arguments.\textquotedblright\  Franklin (1978) urges agents to list
\textquotedblleft strengths and weaknesses\textquotedblright\ to help
decision-making. In political science, \textquotedblleft muddling through"
behaviors (Lindblom, 1959) describe administrative and political
decision-making processes. Agents and organizations are never supposed to
optimize.

Incremental behaviors are state-dependent when choices are sequential,
restricting the choice set by successive eliminations. In most cases,
choosing to do something is anchored in previous choices because choosing to
change entails choosing not to stay. Choosing is rarely between new
options, but between a new and an old one (an anchoring effect of successive
comparisons): to stay or to move, to stop doing something or to carry on, to
buy the same good as before or not. Most behaviors are reference-dependent,
including anchoring, mental accounting, and bracketing.

Repeated choice matters as a step in a dynamic process. This leads to the
concept of temporary routines where the agent repeats the same choice and
changes from time to time, moving from a temporary routine to a new one, to
finally end in a permanent routine.

At one extreme there are incremental behaviors, at the other 
traditional global optimization. In this
ideal case, the agent discovers the whole state space of alternatives in one
preliminary hidden step, before, in one shot, pairwise comparing
alternatives. This represents a static and global cost-benefit analysis.

\smallskip

Our purpose is to pave the way between these extremes, from incremental and
low goal-oriented behaviors to increasingly goal-oriented (instrumental)
behaviors, up to optimizing behaviors. We  introduce intermediate
goal setting \textquotedblleft decision-making" processes, such as
intermediate satisficing (improving enough) where goals change along the
adaptative process. They vary with motivation and depend on goal-setting
costs. Intermediate payoffs become endogenous to the course of \
\textquotedblleft decisions followed by actions\textquotedblright\ and cannot be
cumulated ex ante, as this is done in a substantive dynamical model.

Assume that, out of ignorance, the agent does not particularly optimize (he does neither
know the underlying state space nor his utility function ex ante). At each
step, he must explore around the current state to discover his local
environment and possibly improve his performance. This is a local search
optimizing process (hill climbing). If the agent wants to improve his
motivation to try to \textquotedblleft improve more\textquotedblright\ and
\textquotedblleft more quickly\textquotedblright\ at each step, he can set
intermediate goals (intermediate aspiration and intermediate satisficing
levels). These goals help him to drive his intermediate exploration process
to \textquotedblleft explore enough\textquotedblright\ but \textquotedblleft
not too much\textquotedblright\ around to be able to \textquotedblleft
improve enough\textquotedblright\ at each step. This is the
\textquotedblleft gradual satisficing\textquotedblright\ process of
Soubeyran (2006), Martinez-Legaz and Soubeyran (2002), without inertia and
frictions.

In the presence of inertia (\textquotedblleft costs to
move\textquotedblright ), the agent must do more to sustain his intermediate
motivation: he must \textquotedblleft improve even more\textquotedblright\
at each step to be able to compensate intermediate \textquotedblleft costs
to move\textquotedblright\ by intermediate \textquotedblleft advantages to
move\textquotedblright . Then, it is \textquotedblleft worthwhile to
move\textquotedblright. This is a way to limit intermediate sacrifices.

\smallskip

Our modelling of an incremental decision-making involves three interrelated
blocks:

1) a motivation building and goal-setting block where, at each step, the
agent tries to \textquotedblleft improve enough" (improving or satisficing)
his per unit of time payoff $g(y)\in \mathbb{R}$ generated by each action $%
y\in X$. The agent solves, at each step, a qualitative inequality

\begin{equation}  \label{Intro1}
y\in I \left[x,\varepsilon (x)\right] = \left\{ y\in X:\mbox{ }g(y)-g(x)\geq
\varepsilon (x)\geq 0\right\} \subset X
\end{equation}

\noindent where the temporary satisficing gap is $\varepsilon (x)\geq 0.$

2) a learning block, where, at each step, the agent tries to find a
satisfactory action. This requires to evaluate the per unit of time payoff $%
g(y)$ around the present action $x\in X$. At each step, the agent solves the
qualitative inclusion

\begin{equation}  \label{Intro2}
y\in E\left[ x,r(x) \right] \subset X
\end{equation}

\noindent with $r(x)$ \ equal to the radius of the exploration set around $x\subset X.$

3) a transition \textquotedblleft worthwhile-to-move" block where, at each
step, the agent considers to move or to stay: acceptable changes are such
that the estimated behavioral advantages to move $A(x,y)$ are higher than
some fraction $1\geq \xi (x)\geq 0$ of the estimated behavioral costs to
move $C(x,y)\geq 0$. The agent solves a qualitative \textquotedblleft
worthwhile-to-move" inequality

\begin{equation}  \label{Intro3}
y\in W(x)=\left\{ y\in X:\mbox{ }A(x,y)\geq \xi (x)C(x,y)\right\} .
\end{equation}

Each of these blocks contains three aspects:

i) At each step, \textit{heuristics} and qualitative tools are used by the
agent in order to estimate advantages and costs;

ii) At each step, \textit{punctuated} dynamical aspects enter into the
description of the articulation between static phases of
exploration-exploitation and dynamic moving phases linking two consecutive
temporary routines.

iii) \textit{Behavioral} aspects include physical, physiological,
psychological, cognitive, and social features.

\smallskip

Rumelt (1990) provides an excellent literary presentation approach of
inertia. Conlisk (1996) emphasizes the importance of deliberation costs.
Lipman and Wang (2000, 2006) introduce costs to change actions in the theory of
non adaptive repeated games with perfect rationality (agents optimize a long
term objective, a mean of cumulated payoffs). Our formalization is far more
general, using the concept of distance as a dissimilarity index. We model:

\begin{itemize}
	\item Adaptive processes of decision-making such as \textquotedblleft muddling
through" (Lindblom, 1959) and \textquotedblleft satisficing process" (Simon,
1955) which limit intermediate sacrifices necessary to reach the final
moving goal. Such \textquotedblleft worthwhile-to-move" processes adapt more
or less to inertia during the transition from the beginning to the end. To
our knowledge, there exists no formal model of \textquotedblleft muddling
through."  For a dynamic model of satisficing, see Selten (1998); for a
static model of satisficing, see Tyson (2005).

	\item  Qualitative heuristics of search and exploration, in a topological context,
using an enclosing principle (using inequalities rather than equalities,
putting bounds on control variables) (Gigerenzer and Todd, 1999; Gilovitch et al., 2002).

	\item Path dependency and lock-in effects. We show that \textquotedblleft muddling
through" behaviors are path dependent because they have reference-dependent
payoffs, and make small steps because they converge, due to high costs of
change and bounded needs. The closest model to ours is the local search hill
climbing algorithm (Aarts and Lenstra, 2003) which at each step improves by
local exploration but which ignores inertia.

	\item  Habits and routines as the outcomes of habituation and routinization of
temporary habits and routines which converge to permanent habits and
routines. Our model shows when \textquotedblleft temporary satisficing"
behaviors converge to a \textquotedblleft permanent" routine (local maximum,
attractor, fixed point). This occurs indeed when the agent is motivated
and \textquotedblleft improves enough" at each step. For optimization of
habits and addiction see Abel (1990) and Carroll (2001), for habits and
procrastination see O'Donoghue and Rabin (1999), at the organizational level of
routines see Nelson and Winter (1997).

\end{itemize}

We shall also calibrate the dynamic inefficiency of such adaptive processes
to know how far from optimization agents behave, show how to overcome
inertia and prove the \textquotedblleft (epsilon)- Variational
principle\textquotedblright\ of Ekeland (Attouch and Soubeyran, 2006). It
leads us to revisit the \textquotedblleft
as-if-hypothesis\textquotedblright\ where, as in economics, agents
optimize,\ and put to the fore the \textquotedblleft local search and
proximal algorithms\textquotedblright. These algorithms, which involve
inertia features, are a mixture of local search algorithms and proximal
algorithms. Moreover, it provides a general framework and a large field of
applications for proximal algorithms (Attouch and Bolte, 2006; Attouch et
al., 2007; Attouch and Teboulle, 2004). Connections with second order dynamic
optimization models with memory can be made (Attouch et al., 2000; Attouch
and Soubeyran, 2006; and, for long memory effects, Goudou and Munier, 2005).

\smallskip

In section 2, we consider incremental worthwhile-to-move behaviors used to
bound intermediate sacrifices. In section 3, we emphazise the inertia
context with intermediate costs to move. A typology of costs to move is
introduced. In section 4, we detail the goal-setting context with
intermediate advantages to move. In section 5, we model the punctuated
exploration-exploitation and moving process of decision-making. In section
6, we examine behavioral dynamics with not too many intermediate sacrifices.
In section 7, we consider incrementalism through enclosing. In section 8, we
state the \textquotedblleft worthwhile-to-move\textquotedblright theorem
which gives general conditions for the convergence of the process toward a
permanent routine. In section 9, we revisit the \textquotedblleft
as-if-hypothesis\textquotedblright \ by introducing inertia aspects and put
to the fore the role of \textquotedblleft local search proximal algorithms.''

\section{Incrementalism: \textquotedblleft
Worthwhile-to-Move\textquotedblright \ Behaviors \newline
Bound Intermediate Sacrifices}

\subsection{The State-Dependent Balance between Intermediate Advantages and
Costs to Move}

We model \textquotedblleft muddling through\textquotedblright \ behaviors
(making small steps, improving step by step, Lindblom, 1959) comparing
advantages and costs to change at each step.

\smallskip

The simplest way to decide what to do between two possibilities is to
compare advantages and costs. The point is that most decisions are
state-dependent, part of a dynamic decision-making process. They represent
intermediate decisions to reach well defined goals.

Consider the agent at $x\in X$ who has some motivation to change. He wants
to move from $x\in X$ to some $y\in X$. He explores around $x\in X$, in an
exploration set $E\left(x,r(x)\right)$ of size $r(x)\geq 0,$ to estimate and
compare state-dependent intermediate advantages $A(x,y)\in\mathbb{R}$ and
costs $C(x,y)\in \mathbb{R}^{+}$ to move from $x$ to $y$. The size of an 
exploration set $E(x)$ can be its radius $r(x) = \sup\left\{d(x,y), \  y\in E(x)\right\}.$ 
To simplify, we
consider only scalar advantages and costs. They are reference-dependent, the
reference being, at each step, the state $x$ of departure, then $y$, and so
on. The case of multidimensional advantages and costs to move $A(x,y)\in V$, 
$C(x,y)\in V$ was examined in Soubeyran and Soubeyran (2004 in a Riesz space 
$V$, revised version 2007 in a semi-group $V$).

According to the \textquotedblleft worthwhile-to-move'' principle, an
acceptable move is such that the estimated advantages are higher than some
proportion $1\geq \xi (x)\geq 0$ of his estimated costs:

\begin{equation}  \label{WTM}
A(x,y)\geq \xi (x)C(x,y).
\end{equation}

\noindent A move satisfies the \textquotedblleft worthwhile-to-move''
principle if and only if it satisfies this inequality. The sacrificing rate
is $1-\xi (x)$, the portion of the costs to move which the agent does not
put in the balance is $\left(1-\xi (x)\right)C(x,y)$. The non sacrificing
rate is $\xi (x)$.

The \textquotedblleft \textit{worthwhile-to-move}'' set at $x\in X $ is
defined by

\begin{equation}  \label{WTMset1}
W(x)=\left\{y\in X: \mbox{  } A(x,y)\geq \xi (x)C(x,y)\right\} \subset X.
\end{equation}

This defines a \textquotedblleft worthwhile-to-move'' set-valued
relationship $x\in X\longmapsto W(x)\subset X.$ We assume that
state-dependent advantages and costs to move are zero if the agent stays at $%
x\in X: A(x,x) = C(x,x) = 0$ for all $x\in X.$ Thus $x\in W(x)$ for all $%
x\in X.$ We also assume that $C(x,y)>0$ for $y\neq x.$

\noindent We reformulate this relationship by using a standard goal function 
$g: x\in X\longmapsto g(x)\in\mathbb{R}$. The real number $g(x)$ represents
the instantaneous utility of the agent at state $x$. Starting from a known couple $%
(x,g(x))\in X\times\mathbb{R},$ the agent does not know the values $g(y)$ of
his utility function, for states $y\neq x$, without exploration around $x.$
Intermediate advantages to move are

\begin{equation}  \label{Interad}
A(x,y)= t(y)\left( g(y)-g(x)\right)
\end{equation}

\noindent where $t(y)\geq 0$ is the length of time during which the agent
can benefit or choose to exploit his instantaneous advantages to move $%
g(y)-g(x)\geq 0$. This length of time $t(y)$ defines his intermediate
exploitation. Because costs to move are non negative a \textquotedblleft
worthwhile-to-move'' choice $y\in W(x)$ improves: $y\in W(x)\Longrightarrow
g(y)\geq g(x).$ Costs to move can be zero, in absence of friction. In this
case of no inertia, if the exploitation time $t(y)>0$ is strictly
positive, the worthwhile-to-move relationship reduces to improving: $y\in
W(x)\Longleftrightarrow g(y)\geq g(x).$

\noindent We consider three kinds of inefficiencies: \textquotedblleft lack
of motivation,''   \textquotedblleft lack of knowledge,''  and \textquotedblleft
high local costs to move (inertia)'' which we will model with the help of
three blocks and their interrelations: \textquotedblleft goal setting,'' 
\textquotedblleft exploitation-exploration,''  and \textquotedblleft moving.''
We first examine the moving block, consisting of acceptable transitions
where sacrifices are not too high.

Intermediate sacrifices are maximal when the non sacrificing ratio $\xi
(x)=0 $ is zero. Then, the agent bears costs to move, but ignores them
during the transition. In this case, if the exploitation time is strictly
positive, the worthwhile-to-move relationship reduces once again to
improving, $y\in W(x)\Longleftrightarrow g(y)\geq g(x)$ from $A(x,y) =
t(y)\left(g(y)-g(x)\right) \geq \xi (x)C(x,y)=0$ and $t(y)>0%
\Longrightarrow g(y)\geq g(x).$

When $0<\xi (x)<1,$ the agent ignores the costs to move, or intermediate
sacrifices $\left(1-\xi (x)\right) C(x,y).$ In this case, the ratio
\textquotedblleft advantages over costs to move'' $A(x,y)/C(x,y)$ must be
higher than $\xi (x)<1.$

When $\xi (x)\geq 1,$ the agent is much more demanding, he refuses
intermediate sacrifices, and requires to make an extra gain $A(x,y)-\xi
(x)C(x,y)\geq 0.$ In this case, the ratio \textquotedblleft advantages over
costs to move'' $A(x,y)/C(x,y)$ must be higher than $\xi (x)>1.$

In each case, the \textquotedblleft worthwhile-to-move'' condition is a
reference state-dependent satisficing condition : $A(x,y)/C(x,y)\geq \xi
(x)\geq 0$ for $y\neq x$ and $C(x,y)>0.$ The \textquotedblleft
worthwhile-to-move'' definition (\ref{WTMset1}) can be equivalently
formulated as

\begin{equation}  \label{WTMset2}
W(x) = \left\{y\in X: \mbox{ } A(x,y)/C(x,y) \geq \xi (x)\geq 0 \mbox{ if }
y\neq x\right\} \cup \left\{ x\right\}.
\end{equation}

To simplify this initial presentation of the \textquotedblleft
worthwhile-to-move'' principle, take the intermediate exploitation time $%
t(y)=T$ as a given constant. We will remove this simplification very soon.

\smallskip

The \textquotedblleft worthwhile-to-move'' principle defines an acceptable
transition process $x_{n+1}\in W(x_{n}),\mbox{ } n\in {\mathbb{N}} $. At
each step, when moving from $x_{n}$ to $x_{n+1}$, intermediate sacrifices
are not too high, advantages to move $A(x_{n},x_{n+1})$ are greater than
some fraction $\xi (x_{n})\geq 0$ of costs to move $C(x_{n},x_{n+1})$. At
each step the agent improves his goal from $g(x_{n})$ to $g(x_{n+1}) >
g(x_{n}), \mbox{ } n\in {\mathbb{N}} .$

If the agent wishes to move toward a final goal, and if costs to move are
non negligible, the agent must, at each step, compensate local costs to move $%
C(x,y)$ by some local advantages to move $A(x,y).$ This defines an
acceptable transition process.

\smallskip

The \textquotedblleft worthwhile-to-move'' principle governs a lot of weakly
goal-oriented behaviors, where, at each step, an agent, starting from the
current state, tries both to improve, and to balance local advantages to
move to a minimal fraction of costs to move. We will show that the course of
choices and actions of this local \textquotedblleft 'behavioral
cost-benefit'' decision-making is localizable, because it is nested. But it
does not necessarily shrink to a point. Such a process is incremental,
reference and path-dependent. This \textquotedblleft worthwhile-to-move''
process characterizes a behavioral \textquotedblleft cost-benefit'' analysis
and a \textquotedblleft muddling through'' process of administrative
behavior (Lindblom, 1959), where frictions and inertia play a major role.
This incremental decision-making describes a\ \textquotedblleft mental
accounting'' behavior.

\smallskip

Consider an agent who, at state $x$, benefits of the instantaneous utility $%
g(x)$ while his aspiration level at $x$ is higher: $\widehat{g}(x) > g(x).$
While staying at $x\in X$, he feels the frustration $f(x)=\mu (x)\left( 
\widehat{g}(x)-g(x)\right) >0$ generated by the per unit of time unsatisfied
need $\widehat{n}(x)=\widehat{g}(x)-g(x)>0.$ The weight $\mu (x)>0$
represents his per unit frustration feeling. This pushes him to change for a
new state $y\in X$ with a better instantaneous utility $g(y)>g(x)$. This is
incrementalism where an agent tries to gradually improve his instantaneous
utility. This is the case of a goal-oriented agent who only wants to reduce
his frustration feeling by improving. He will either explore in a given
exploration set or enlarge his exploration set step by step. In both cases,
he will stop the exploration around $x$, in the exploration set $E\left(
x,r(x)\right) \subset X$ as soon as he finds a new state $y\in E\left(
x,r(x)\right) $ which improves,\ $g(y)>g(x).$

This is the no-friction case with no cost to move. In the friction case, the
agent must compensate the costs to move by high enough advantages to move.
He must improve enough to limit the intermediate sacrifices of moving. This
defines a \textquotedblleft worthwhile-to-move'' constraint during the
transition, $y\in W(x),$ which is expressed as

\begin{equation}  \label{WTMC}
T\left( g(y)-g(x)\right) \geq \theta (x,y)d(x,y)
\end{equation}

\noindent where $d$ \ is a distance on $X$, the intermediate exploitation time 
$t(y)=T$ is taken as a given constant. In our simple case, we take the
ratio $\theta (x,y)$ higher than some constant $\theta >0.$

As said before, incrementalism in terms of goal improvement is when an agent
tries to improve his per unit of time utility step by step. A world of
frictions adds \textquotedblleft costs to move.''  There is a need to
\textquotedblleft improve enough'' at each step in order to limit
intermediate sacrifices coming from costs to move. We will show that such
\textquotedblleft improving enough'' process generates local actions, small
steps, local \textquotedblleft worthwhile-to-move'' transitions toward a new
improving state. These improving transitions must be acceptable, in the
sense that they present not \textquotedblleft too many sacrifices'' in the short
run. In this case, exploration is not goal-oriented. The agent learns by
doing.

A worthwhile-to-move behavior usually requires more than improving $%
g(y)-g(x)>0$. This is \textquotedblleft improving enough'', $g(y)-g(x)\geq
\varepsilon (x)>0$, for some given and feasible gap of improvement $%
\varepsilon (x)>0,$ to compensate for moving costs. This requires both $%
T\left( g(y)-g(x)\right) \geq \varepsilon (x)\geq 0$ and $\varepsilon
(x)\geq \theta (x,y)d(x,y).$ This remark will help us later to make the link
between a worthwhile-to-move process and an intermediate satisficing process.

\subsection{Local Actions and Convergence}

We  show that an incremental \textquotedblleft worthwhile-to-move''
behavior has the \textquotedblleft local action'' property, and converges
(both goals and states). Local action is a consequence, not an hypothesis.

The state space $X$ is a metric space with distance $d.$ Take the
exploitation time $T=1$, the \textquotedblleft worthwhile-to-move''
inclusion (\ref{WTMC}) becomes 
\begin{equation*}
W(x)=\left\{ y\in X: \mbox{  } g(y)-g(x) \geq \theta (x,y)d(x,y)\right\}
\subset X.
\end{equation*}

We suppose that the agent has enclosed the \textquotedblleft
worthwhile-to-move'' inclusion $y\in W(x)$ within the inclusion $%
S(x)\supseteq W(x)$ in the following way. The agent can control the process
step by step, in such a way that, for all $x,y\in X,$ $\theta (x,y)\geq
\theta >0.$ Let

\begin{center}
$S(x)=\left\{ y\in X: \mbox{  }g(y)-g(x)\geq \theta d(x,y)\right\}. $
\end{center}

With $x_{0}\in X$ being given, the \textquotedblleft worthwhile-to-move''
process $x_{n}\rightarrow W(x_{n}), \mbox{ } n\in {\mathbb{N}}$, is enclosed
in the intermediate satisficing process $x_{n}\rightarrow S(x_{n}), \mbox{ }
n\in {\mathbb{N}}$.

The relationship $S$ defined by $xSy \Leftrightarrow y\in S(x)$ is reflexive
and transitive:

i) $x\in S(x)$ for all $x\in X$ is a consequence of \ $d(x,x)=0$ \ for all $x\in
X.$

ii) $y\in S(x)$ and $z\in S(y)$ implies $z\in S(x)$. This is a consequence
of the triangle inequality: $d(x,z)\leq d(x,y)+ d(y,z)$.

Indeed, from $g(y)-g(x)\geq \theta d(x,y)$ and $g(z)-g(y)\geq \theta d(y,z)$%
, by adding the two inequalities, $g(z)-g(x)\geq \theta d(x,z)$.

Hence the enclosing process is nested: $S(x_{0})\supset S(x_{1})\supset
....\supset S(x_{n})\supset S(x_{n+1})\supset ...$

This allows the agent to gradually improve the localization of the
\textquotedblleft worthwhile-to-move'' process: Let $\left\{x_{n+1}\in
W(x_{n}), \mbox{  } n\in \mathbb{N} \right\}$ with $x_{0}\in X $ be a
\textquotedblleft worthwhile-to-move'' process which can be enclosed, using
the enclosing heuristic

\begin{center}
$W(x)\subset S(x)=\left\{ y\in X: \mbox{  }g(y)-g(x)\geq \theta
d(x,y)\right\} $ for all $x\in X$.
\end{center}

Then the \textquotedblleft worthwhile-to-move'' process is nested.

\smallskip

The agent has limited needs, his utility function $g: x\in X\mapsto
g(x)\in\mathbb{R}$ is upper bounded, set 
$\mbox{  } \overline{g} = \sup_{x\in X}  g(x) <+\infty $. 
Assume that the costs to move are high enough
locally: $C(x,y)\geq \theta d(x,y),\mbox{  }\theta >0$. Then, during the
enclosing process $x_{n}\rightarrow S(x_{n}),\mbox{  }n\in {\mathbb{N}} $
with $x_{0}\in X $ given, we have the succession of inequalities

\begin{equation}  \label{WTMn}
g(x_{n+1})-g(x_{n})\geq \theta d(x_{n},x_{n+1}),\mbox{  }n\in {\mathbb{N}}.
\end{equation}

Let us state the simple but fundamental local action property:

\smallskip \textbf{The Local Action Proposition: } \textit{If the agent has limited
needs, high local costs to move, and if he accepts limited intermediate
sacrifices, by enclosing, this implies local actions and convergent goals}:

\begin{center}
$d(x_{n},x_{n+1})\rightarrow 0$ \ \textit{as} \ $n\rightarrow +\infty $ \ \textit{and} \ 
$g(x_{n})\rightarrow g^{\ast }\leq \overline{g}<+\infty .$
\end{center}

\textit{Proof}: By inequality (\ref{WTMn}), the sequence of goals $\left\{
g(x_{n}),\mbox{  }n\in {\mathbb{N}}\right\} $ is increasing (by definition
of a distance, $d(x_{n},x_{n+1})$ is nonnegative). As needs are limited, ($%
\overline{g}<+\infty $), this sequence converges to some limit $g^{\ast
}\leq \overline{g}<+\infty $. It follows from Eq.(\ref{WTMn}) that the
distance between two successive states tends to zero, $d(x_{n},x_{n+1})%
\rightarrow 0$ as $n\rightarrow +\infty $.

\smallskip Assuming that the state space is complete, we obtain the
convergence property: \smallskip

\textbf{The Convergence Proposition: }\textit{If the state space $X$ is a complete
metric space, if the per unit of time utility function $g(.)$ is upper
bounded, if the worthwile to move process $x_{n}\rightarrow W(x_{n}), 
\mbox{
}n\in{\mathbb{N}} $ is enclosed within the enclosing process $%
x_{n}\rightarrow S(x_{n}),\mbox{  }n\in{\mathbb{N}}$ with $x_{0}\in X$ given
and high enough local costs to move, then, the \textquotedblleft
worthwhile-to-move'' process converges toward some final state, $%
x_{n}\rightarrow x^{\ast }$ as $n\rightarrow +\infty $.}

\smallskip

\textit{Proof}: By adding all the inequalities (\ref{WTMn}) from $n=0$ to $%
m, $ using the subadditivity of the distance function, and the fact that the
instantaneous utility function is upper-bounded  ($\overline{g}<+\infty $),
 one obtains
\begin{equation}  \label{WTMsum}
+\infty  > \overline{g}-g(x_{0}) \geq  g(x_{m+1})-g(x_{0}) \geq  \theta
\sum_{n=0}^{m}d(x_{n},x_{n+1}).
\end{equation}

As a consequence, the series \  $\sum_{n=0}^{+\infty }d(x_{n},x_{n+1})<+\infty $
is convergent. The \textquotedblleft worthwhile-to-move'' sequence
 $\left(x_{n}\right)_{n \in \mathbb{N}} $\ is a Cauchy sequence which converges to some state $x^{\ast
} $ $\in X$ in the complete metric space $X.$

\subsection{From Rest Points to Behavioral Rest Points}

A performance $x^{\ast }\in X$ is said to be a \textit{rest point} if, for
any $y\in X$ with $y\neq x^{\ast }$, it is not worthwhile to move from $%
x^{\ast }$ to $y$. This is equivalent to say that $x^{\ast }\in X$ is a rest
point of the \textquotedblleft worthwhile-to-move''\ relationship $x\in X$ $%
\longmapsto W(x)\subset X,$ $W(x^{\ast })=\left\{ x^{\ast }\right\}.$ In
this case, the agent has no further incentive to move. Thus, $x^{\ast }\in X$
is a rest point iff for any $y\neq x^{\ast },$ $A(x^{\ast },y)<\xi (x^{\ast
})C(x^{\ast },y)$. In our simplified example, intermediate advantages to
move are proportional to per unit of time advantages to move, $A(x,y)=
T\left( g(y)-g(x)\right) =g(y)-g(x). $ If the exploitation time is given, $%
T=1$, and costs to move are proportional to the distance of moving, $%
C(x,y)=\theta d(x,y)\geq 0, \mbox{  }\theta >0$, then $x^{\ast }\in X$ is a
rest point if, for any $y\neq x^{\ast } $,

\begin{equation*}
g(y)-g(x^{\ast })<\theta d(x^{\ast },y).
\end{equation*}

This property defines a rest point in a strong sense.

A weak rest point $x^{\ast }\in X$ is such that $g(y)-g(x^{\ast })\leq
\theta d(x^{\ast },y),$ for all $y\in X$, with $\theta \geq 0.$ If $\theta
=0,$ a weak rest point is a global maximum, while a strong rest point is the
unique global maximum of the function $g(.).$ The traditional optimization
problem to solve $\sup_ {x\in X}  g(x)$ is a particular case.

A performance $x^{\ast }\in X$ is said to be a \textquotedblleft \textit{%
behavioral rest point}" with respect to a given initial data $x_{0}\in X,$
if it is a rest point, and starting from the given state $x_{0}\in X,$ it
can be reached, following a \ \textquotedblleft worthwhile-to-move''
(acceptable with not too many intermediate sacrifices) transition process $%
x_{n}\rightarrow W(x_{n}),\mbox{  }n\in{\mathbb{N}}.$

Substantive rationality considers only the case of a rest point $x^{\ast
}\in X$, more precisely a global or a local maximum. If the agent happens to
be there, he prefers not to deviate.

The Ky Fan theorem (Fan, 1972) gives conditions over the net incremental
gain \textquotedblleft to move instead of to stay" $\Delta (x,y) =
A(x,y)-\xi (x)C(x,y)$ which guarantee the existence of a rest point. The Ky
Fan conditions are (Singh et al., 1997: 137): Let $X\subset \Upsilon $ be a
non empty convex set in a topological vector space $\Upsilon $. Let $\Delta
:(x,y)\in X\times X\longmapsto \Delta (x,y)\in $ be such that:

\begin{enumerate}
\item for each $x\in X,$ $\Delta (x,y)$ is a quasi-concave function of $y\in
X,$

\item for each $y\in X,$ $\Delta (x,y)$ is a lower semicontinuous function
of $x\in X,$

\item $\Delta (x,x)\leq 0$ for all $x\in X$,

\item $X$ is compact.
\end{enumerate}

\noindent Then there exists a point $x^{\ast }\in X$ such that $\Delta
(x^{\ast },y)\leq 0$ for all $y\in X.$

But this story does not tell us why the agent is there, how he has reached
such a rest point. Procedural rationality considers the much more convincing
case of a behavioral rest point. The story tells us where the agent starts
from, and, step by step, which acceptable paths he is supposed to follow to
reach or to be locked in a behavioral rest point. Procedural rationality
examines which kind of process converges in a rest point, starting from any
initial position $x_{0}\in X$, do rest points exist, and does such a process
converge in finite time? In the next sections, we show that, under some mild
conditions, the \textquotedblleft worthwhile-to-move'' process \ $x\in
X\longmapsto W(x)=\left\{ y\in X, \mbox{ } \Delta (x,y)\geq 0\right\}
\subset X$ converges to a rest point. Convergence in finite time requires
more.

Concrete examples of behavioral rest points are daily life behaviors such as
routines, habits, rules, norms, practices, lock-in effects.
A behavioral rest point can be inefficient because it is constrained by a
transition path which ends to it. How far from the supremum does an
improving process lead? The inefficiency gap of such a behavior will be
measured by the difference $\overline{g}-g^{\ast }\geq 0$ between the
supremum $\overline{g}<+\infty $ of \ $g(.)$ over the state space $X$ and
the limit $g^{\ast}=\lim_{n\rightarrow +\infty }g(x_{n})\leq \overline{g} $
of the improving process.

\section{The Inertia Context: \textquotedblleft Intermediate Costs to
Move\textquotedblright}

We have first examined the balance between behavioral advantages and costs
to move, to explain how incremental behaviors (muddling through processes)
emerge (converge, hence making small steps). Before examining more
goal-oriented behaviors, we examine the two terms of the balance, behavioral
advantages and costs to move (which are not the usual benefit and cost
functions) including goal-setting, cognitive, psychological, and inertial
aspects. In this section we first examine costs to move.

\subsection{Frictions and Inertia}

Management sciences consider distorted perception, dulled motivation,
failed creative response, political deadlocks, and action disconnecting
(Rumelt, 1990). The term inertia is used for each of these frictions.

\begin{itemize}
\item Distorted perception comes from myopia, (the unability to forecast the
future with clarity), and from denial (a defensive behavior, which is the
rejection of information which is contrary to what is desired or what is
believed to be true). Denial may stem from hubris or from fear. Information
filtering rejects information which is unpopular, unpleasant, or contrary to
doctrine. Grooved thinking rejects information which deviates too much from
common wisdom.

\item Dulled motivation describes the lack of sufficient motivation to
change, due to the abandonment of costly sunk specific investments.

\item Failed creative response concerns the difficulty to choose a direction
because of the complexity of the choice, the speed of change, inhibition and
non-reactive state of mind, or inadequate strategic vision.

\item Political deadlocks come from the three main sources of disagreement
among agents: difference in personal interest, in beliefs, and in values.

\item Action disconnecting comes from lack of vision,
lack of leadership, attachment to statu quo, embedding in routines where
changing one part of the process requires to change a lot of other parts.
\end{itemize}

To capture most of these frictions we define

\begin{itemize}
\item i) Motivation building, goal-setting, and exploitation (doing,
know-how).

\item ii) Exploration (information seeking, concept acquisition).

\item iii) Moving (changing, learning).
\end{itemize}

These three kinds of frictions introduce reference-dependent effects on
goals and costs. The terminology of inertia is limited to the moving block,
to costs to move more or less quickly.

\subsection{Costs to move and costs to stay}

Costs to move can be either physical or behavioral. Costs to stay and dynamical
costs to move can be classified into: costs linked to a task (costs to
start, to progress, to stop), costs to switch between tasks, and costs of
interactions between agents.

\smallskip

1. Costs to Improve the Way of Doing a Task (Intra Costs):

\begin{itemize}
\item Motivation costs are costs to set a goal or strive for it, to set an aspiration
(longing) level, a satisficing level, and to adjust the longing level.
Excitation costs are often necessary to start some novel action. Emotions
drive the goal setting process.

\item Physical and physiological dynamic costs
can be fixed costs to start an action, such as training, preliminary, or
warming-up costs to repeat a task, or costs to initiate a new way of doing a
task. There are also costs to stop an action. Dynamic costs to change concern
inertia and reactivity. They increase with the speed with which the
agent wants to do something, or to do it more quickly. They represent
dynamic production costs which are linked to the intensity of effort and the
speed of doing. 
Reactivity is one of the most important feature of modern organizations.

\item Costs to stay, maintenance and operational costs are costs to repeat
an action. They are related to the difficulty of the way of doing a task,
and also to boredom and displeasure. Weariness appears with repetitions and
successions of several actions, without enough time to recover.

\item Cognitive costs concern knowledge acquisition, such as exploration
costs or costs to acquire information and knowledge. Cognitive costs concern
deliberation costs, costs to choose, to eliminate, to renounce, costs of
knowledge acquisition. Knowledge costs are costs to solve a problem. They
depend on the residual difficulty to solve this problem. Deliberation costs
before action belong to the category of exploration costs (Conlisk, 1996).
Learning costs can be behavioral costs to learn how to do  a task differently
in a better way, they are costs of \textquotedblleft ways of doing'' like
recipe costs, and costs to \textquotedblleft learn how to learn''. Costs of
sampling and gathering information by trials and errors are prominent at the
beginning of a learning process. The agent must determine if a choice is an
error or not.

\item Psychological costs include costs to direct attention, costs coming
from emotions, stress, anxiety, costs of  doubt,
ambiguity and risk, costs to accept to choose early, to eliminate possible
solutions, or to renounce to some opportunity, to fear to regret without
enough knowledge and information.

\item Intertemporal  or transition costs are linked to the costs of not
reaching a target immediately. They represent the costs of intermediate
sacrifices which the agent has to accept before reaching a final goal. These
costs include impatience, resilience, boredom from repetition, stress, costs
to change to often, costs to wait, to delay, costs to accept to be mistaken.
\end{itemize}

Costs to move take very different shapes, depending on the nature of the
task: consuming, working, coaching, managing.

\smallskip

2. Costs to Switch between Tasks (Inter Costs):

Costs to switch from one task to another concern costs of dissimilarity,
and, as before, costs to start actions, and costs to stop (inhibition
costs). Costs of doing new actions increase with the degree of dissimilarity
with past actions. Ehrlich (1975) insists on the fact that similar actions
following a given action are much easier to do than disconnected actions
belonging to different fields of competence. Costs of doing several tasks at
the same time and costs of specialization (costs of doing the same task
again and again) are intermingled. Inhibition costs are the costs to forget
and leave a given task to be able to start a new task. Inertia costs can be
addiction costs to stop consuming drugs, and costs to escape from habits. The
more often an agent changes, the easier he changes again, or the more
tiresome it is.

Adaptation and adjustments costs regroup costs to switch from one task to
another, when an agent has to make several tasks (extensive aspect) and
costs to change the \textquotedblleft way of doing'' a given task (intensive
aspect) if he makes only one task. These costs to change the
\textquotedblleft way of doing'' a task are behavioral.

Switching costs (Klemperer, 1995) concern adaptation costs for consumer
tasks. For a consumer, to consume a good is a task, a kind of production
act, a consumption recipe. Switching costs are adaptation costs for a
consumer who chooses to move from the consumption of one good to another.
They concern costs to lose compatibility with his existing physical
equipment, his current knowledge, relational costs to switch from a supplier
to a new one, costs to lose discount coupons, costs to learn using new
brands, psychological costs to try new products.

\smallskip

3. Interaction and Network Social Costs:

Costs to stay and costs to change concern also interaction costs between
agents, notably, network costs to keep, build, or delete links. 
Transaction costs are costs of making an economic
exchange, such as the effort  required to find out the prefered variety of
a product, the travel time between home and store,
queuing time. Transaction costs include search and information costs,
bargaining costs and policing, monitoring and enforcement costs. Following
Coase (1937) and then Williamson (1975), the determining factors of
transaction costs are frequency, specificity, uncertainty, bounded
rationality, and opportunistic behavior. Transaction costs, or frictions are
dynamic (Langlois, 2005). They concern coordination and synchronization
costs, costs to improve fits, cooperation costs, costs to exchange and build
common knowledge, and organizing costs.

\subsection{Modelling behavioral costs to move and costs to stay}

For an agent having reached a given state $x\in X,$ we model

\begin{enumerate}
\item Costs to stay or maintenance costs (Operational costs):

If the agent stays at $x$, he must pay maintenance costs $m(x)$, or costs to
stay at $x$. Costs to stay at $x\in X$ include the physical costs to repeat $%
x$ and the psychological costs to feel a persistent frustration coming from
the unsatisfied needs $f(x)=\mu (x)\left( \widehat{g}(x)-g(x)\right) >0.$ We
will embed maintenance costs in the instantaneous utility gain, $g(x)=\varphi (x)-m(x)$,
where $\varphi (x)$ is the gross utility.

\item Physical and Behavioral Costs to Move:

If the agent moves, he must pay physical, physiological, psychological, and
cognitive costs of moving from $x$ to some better state $y.$ Let $C(x,y)\geq
0$ be the cost to move from $x$ to $y.$

\item Cognitive Exploration Costs :

If the agent moves, he must spend exploration expenditures $K(x)$ to know
if, after exploration, it will be worthwhile to move from $x$ to some better
state $y.$ At each step, and before moving, exploration concerns the estimation of the two
terms of the balance between advantages to move $A(x,y)$ and costs to move $%
C(x,y)$. Exploration costs are static,  moving costs are dynamic.

\item Opportunity Costs to Move:

Costs to move include opportunity costs not to exploit  current
opportunities. If, moving from $x$ to $y$, the agent fails to exploit
his utility $g(x)$ during a time $t(x,y),$  the corresponding opportunity
costs are $O(x,y) = t(x,y)g(x)$.
\end{enumerate}

Because of space limitation, and the huge variety of costs to change, we
give a reduced mathematical form for costs of moving. The set of
alternatives $X$ is a complete metric space. The distance between two
alternatives $x\in X$ and $y\in X$ is an index of dissimilarity between
them. The physical costs of moving from $x$ to $y$ is decomposed into

\begin{center}
$C(x,y) = e(x,y)d(x,y) = t(x,y)c(x,y).$
\end{center}

Costs to move are zero if the agent does not move: $C(x,x)=0,$ or $t(x,x)=0,$
and $c(x,x)=0.$ These formula define the per unit of distance cost to move $%
e(x,y)\geq 0$ and the per unit of time cost to move $c(x,y)\geq 0$. Efforts
per unit of distance are not symmetrical: $e(y,x)\neq e(x,y)$,  \ so costs to
move are not symmetrical either. We can have $C(x,y)\neq C(y,x).$ This
justifies the use of relative entropy such as Kullback-Liebler or Bregman
distances.

With regard to the dependence of the costs to move with respect to the
distance, we distinguish two classes:

\begin{itemize}
\item \textit{High local costs to move} which correspond to a minimum per
unit of distance effort of moving $e(x,y)\geq \underline{e}>0.$ As an
example $C(x,y)= \underline{e} d(x,y)$, whose corresponding mechanical
notion is dry friction. This situation is examined in section 8.

\item \textit{Low local costs to move} which is the complementary class.
This includes the important case $C(x,y)= \underline{e} d(x,y)^2$, whose
corresponding mechanical notion is viscous friction. Because of the presence
of the square, small changes induce very small costs (such as passing from $%
1/10 $ to $1/100$). In that case, convergence is more difficult to prove, it
requires extra geometrical assumptions on the gain or utility function $g$
(such as quasi-convexity or analyticity). A result illustrating this situation
(local search and proximal algorithm) is given in section 9.
\end{itemize}

The dependence of costs to move with respect to time is also a rich topic.
The term $t(x,y)$ is the time spent to move physically from $x$ to $y.$
Among them, examples are reactivity costs to move $C(x,y)=t(x,y)c(x,y)$,
where the instantaneous cost to move $c(x,y)=\rho \left( x,d(x,y)/t(x,y)\right) $ depends
on the mean speed of moving $v(x,y)=d(x,y)/t(x,y)$ (Attouch and Soubeyran,
2006). The cost to move $c(x,y)=\rho \left( x,d(x,y)/t(x,y)\right) $
increases more or less with speed.

\section{The Goal-Setting Context: Intermediate Advantages to Move}

Intermediate advantages to move include intermediate aspiration levels with
psychological aspects. For each state $(x,g(x))\in X \times \mathbb{R}$, \ let 
$\widehat{g}(x)> g(x)$ be the intermediate aspiration level of the agent, which
is an estimation of the unknown supremum $\overline{g}$ \ for example, or
any other level.

Reference-dependent advantages to move: Assume that, after exploration
around the current state $x$ and given his aspiration level $\widehat{g}(x)$,
the agent can find an improving state $y\in E\left( x,r(x)\right) ,$ but
failed to reach his aspiration level, $g(x)<g(\eta)<\widehat{g}(x)$ for all $%
\eta\in E\left( x,r(x)\right).$ At $y\in X,$ the agent will have two
opposite feelings. The gap $\widehat{g}(x)-g(x)>0,$ is a half full or half
empty bottle:

\begin{itemize}
\item He has first a satisfaction feeling $\omega (x,y)=\lambda (x)\left(
g(y)-g(x)\right) $ to have improved, $\ g(y)-g(x)>0$, where $\lambda(x) > 0$
is the weight the agent puts on satisfaction;

\item He has also a disappointment feeling $f(x,y) = \mu (x)n(x,y)$, where $%
\mu(x) > 0$ is the weight the agent puts on disappointment, to have been
unable to satisfy his initial ambition, which decreases with $n(x,y)=%
\widehat{g}(x)-g(y)>0;$

\item The agent can put a weight $\nu g(y)$ on the improving utility $g(y),$
where $\nu =\nu (x)>0;$
\end{itemize}

Instantaneous weighted advantages to move are 
\begin{eqnarray*}
a(x,y)&=&\nu g(y)+\lambda \left( g(y)-g(x)\right) -\mu \left( \widehat{g}
(x)-g(y)\right) \\
&=&\nu g(x)+(\lambda +\nu )\left( g(y)-g(x)\right) -\mu \left( ( \widehat{g}%
(x)-g(x))-(g(y)-g(x))\right) \\
&=&(\lambda +\mu +\nu )\left( g(y)-g(x)\right) +\nu g(x)-\mu \left( \widehat{%
g}(x)-g(x)\right).
\end{eqnarray*}
At each step, the reference is the advantage or disadvantage to stay
\begin{center}
$a(x,x)=\nu g(x)-\mu \left( \widehat{g}(x)-g(x)\right)$
\end{center}
\noindent which is the frustration to stay at $x\in X$. The agent considers
the difference between his advantages to move and his advantages to stay
\begin{equation}  \label{Adv1}
a(x,y)-a(x,x)=(\lambda +\mu +\nu )\left( g(y)-g(x)\right) .
\end{equation}
\noindent This is a reference-dependent criterion for the advantage to move.

\smallskip

The agent may also consider the dynamic advantage to move, which is the
advantage to move multiplied by the duration of hoping to benefit from this
advantage to move. Let $t(y)>0$ be the length of time the agent chooses,
starting from $x$, to stay at $y$ and to exploit the benefit of this new
state $y.$ The temporary \textquotedblleft reference-dependent'' advantage
to move with respect to stay is
\begin{eqnarray}  \label{Adv2}
A(x,y) &=& t(y)\left( a(x,y)-a(x,x)\right) \\
&=& t(y)\delta (x)\left( g(y)-g(x)\right)>0
\end{eqnarray}
where 
\begin{equation}  \label{Adv3}
\delta(x) = \lambda(x) + \mu(x) + \nu(x) > 0
\end{equation}
\noindent is the \textit{character index} of the agent. This formulation
includes  \textquotedblleft improving'', when the agent is not much
goal-oriented, with $\lambda (x)=\mu (x)=0, \mbox{ }\nu (x)>0,$  and
\textquotedblleft improving enough,'' when the agent is much more
goal-oriented, $\lambda (x), \mu (x), \nu (x)>0$.

\section{Punctuated $\mbox { }$\textquotedblleft Exploitation-Exploration''
and Moving Processes}

We define the main ingredients of a generic exploitation-exploration
phase $\left[x\right]$. This notation includes the state $x$ and the control
variables $h(x)$, $\alpha (x)$, $r(x)$ which are described below.

Starting from $x$, the duration $h(x)\geq 0$ of an exploitation-exploration
period is chosen first. Then, the agent chooses how to share $h(x)$ between
exploitation and exploration. At each instant he spends $t(x)=\alpha (x)h(x) 
$ units of time to exploit his instantaneous utility $g(x)$ and $\left( 1-\alpha
(x)\right) h(x)=\tau (x)$ units of time to explore around. Thus, $\alpha
(x)\in \left( 0,1\right) $ and $h(x)\geq 0$ are control variables.

During the phase of exploitation-exploration, the agent stays at $x\in X.$ The longer the
agent exploits the utility $g(x)$, the longer he can explore around, for a
given share of time $1-\alpha (x)$. If the agent, after exploration and
moving, has reached the improved state $x$ with respect to a previous state $%
x_{0}$ such that $g(x)>g(x_{0}),$ and chooses to stay at $x$ during $%
t(x)\geq 0,$ the longer he spends to choose $t(x),$ the larger his
advantages to move from $x_{0}$ to $x,$ $A(x_{0},x)=t(x)\delta (x)\left(
g(x)-g(x_{0})\right) .$ A fixed horizon limits this possibility.

At each step, the agent has some exploration expenditures $K=K(x)=\tau
(x)k(x)\geq 0$ to explore around the current state $x.$ This amount of
resource $K(x)$ helps him discover local advantages $A(x,y)$ and costs to
move $C(x,y)$ around the current state $x.$ As before, to simplify, we take
the exploration expenditure equal to one, $k(x)=1.$ The size of the
exploration set increases with exploration expenditures, $r(x)=r\left(
x,K(x)\right).$

The necessity to repeat exploration at each step comes from the anchoring
aspect of advantages and costs to move$.$ The agent discovers all the
advantages $A(x,y)$ and costs to move $C(x,y)$ from $x$ to $y\in E\left(
x,r(x)\right) ,$ and then moves from $x$ to $y;$ he does not know yet the
advantages and costs to move $A(y,z)$ and $C(y,z)$ from $y$ to $z$. This
forbids him to consider cumulative payoffs as well.

\smallskip

At the end of the exploitation-exploration period $\left[ x\right] $ of
duration $t(x),$ the agent estimates the advantages and costs to move for
all states within the local exploration set: $y\in E\left( x,r(x)\right)
\subset X.$ He compares them, and chooses to stay at $x$ or move from $x$ to
some $y\in E\left( x,r(x)\right).$ If for some $y\in E\left( x,r(x)\right)$,
advantages to move are greater than costs to move (including opportunity
costs to move $O(x,y)$), say $A(x,y)\geq \xi \left( C(x,y)+O(x,y)\right) $,
where $\xi >0$ is a given non sacrificing ratio, he will rather move than to
stay at $x.$ Exploration around is a way for him to be rational. If he takes
the decision to move from $x$ to $y$, the agent will leave the static
temporary routine phase $\left[ x\right] $ to enter into a  phase of moving,
denoted by $\left[ x,y\right]$, which links the previous static
exploitation-exploration phase $\left[ x\right] $ to a new
exploitation-exploration phase $\left[ y\right] .$ A new static period $%
\left[ y\right] $ of exploitation-exploration follows a  phase $\left[
x,y\right] $ of moving.

Given the distance $d(x,y)$ between $x$ and $y,$ the control variables of
the moving phase are its duration $t(x,y),$ the instantaneous cost of moving 
$c(x,y)$, and the per unit of distance cost $e(x,y)$. Costs of moving
\begin{equation}  \label{Cvar1}
C(x,y)=t(x,y)c(x,y)=e(x,y)d(x,y)
\end{equation}
\noindent and the modulus of the speed of moving
\begin{equation}  \label{Cvar2}
\rho (x,y)= \frac{d(x,y)}{t(x,y)}
\end{equation}
\noindent link these choice variables:
\begin{equation}  \label{Cvar3}
c(x,y)=e(x,y)\rho (x,y).
\end{equation}
The context may lead us to add \textit{opportunity costs of moving} $%
O(x,y)=t(x,y)g(x)$, when we consider a single agent, who has a strict time
constraint. These opportunity costs are the lost gains, because the agent is
compelled to stop exploitation during the moving phase of duration $t(x,y).$
For a group, such opportunity costs do not exist because the group can
escape from time constraints and exploit the amount $t(x,y)g(x)$ during this
period of moving. The group can  hire workers to move from $x$ to $y.$ For
a single agent  the balance is
\begin{equation}  \label{Opport}
A(x,y)\geq \xi \left( C(x,y)+O(x,y)\right)
\end{equation}
and $A(x,y)\geq \xi C(x,y)$ for a group.

\section{Behavioral Dynamics: Not Too Many Intermediate Sacrifices}

The \textquotedblleft worthwhile-to-move'' \textquotedblleft decision +
action'' process follows a punctuated dynamic, an alternation of static and
dynamic periods of exploitation-exploration and  periods of moving. The agent
faces contradictory choices, which lead him to use heuristics of choice
before taking action.

\subsection{A Punctuated Dynamic}

Starting from $x$ the agent first chooses the duration $h(x)$ of the
exploitation-exploration period. Then, the agent makes a second choice. At
each time, a fraction of time $\alpha (x)\geq 0$ is devoted to exploitation,
and the other fraction $1-\alpha (x)\geq 0$ to exploration. He exploits
during the duration $t(x)=\alpha (x)h(x) $ and explores during the duration $%
\left( 1-\alpha (x)\right) h(x)$. Consider a punctuated dynamic, a
\textquotedblleft stop and go'' sequence made of three periods starting with
a static exploitation-exploration period $\left[ x\right] ,$ followed by a
 period $\left[ x,y\right] $ of moving and a new static exploitation-exploration
period $\left[ y\right] $ of respective lengths $h(x)\geq 0,$ $t(x,y)\geq 0$
and $h(y)\geq 0.$ The decision to stay at $x$ and to lengthen the period $%
\left[ x\right] $, or to move from $x$ to $y$, entering the moving period $%
\left[ x,y\right] $ to reach the new temporary routine period $\left[ y%
\right] ,$ is, by definition, taken at the end of period $\left[ x\right] .$
At the beginning of the  period $\left[ x,y\right]$ of moving, in order to
decide to move or not from $x$ to $y$, the agent compares the estimated
incremental advantages to move $A(x,y)=t(y)\delta (x)\left( g(y)-g(x)\right)$
 to the estimated costs to move \textquotedblleft without opportunity
costs'' $C(x,y)=t(x,y)c(x,y)$, or costs to move \textquotedblleft with
opportunity costs'' $C(x,y)+O(x,y).$ Exploration costs $K(x)$ are not
involved in the costs of moving $C(x,y)$, but in the larger category of
costs to change.

\subsection{The \textquotedblleft Worthwhile-to-Move\textquotedblright\
Inclusion}

For the sake of simplicity, we consider the case of an agent with
\textquotedblleft no opportunity costs'' to move (the other case will be
examined just after). For an agent following this punctuated dynamic, it is
worthwhile to move from $x$ to $y$ if estimated net advantages to move $%
A(x,y)$ are greater than estimated costs to move $C(x,y)$: $A(x,y)\geq
C(x,y),$ or some portion $0<\xi(x)\leq 1$ of them, $A(x,y)\geq \xi
(x)C(x,y). $ This means that some temporary sacrifices are allowed. In real
life, at each step, the agent accepts to make some temporary sacrifices. He
does not include  some fraction $\left( 1-\xi (x)\right)
C(x,y)$ of the costs to move $C(x,y)$ in the balance. This helps him to hope to improve
more easily from $g(x)$ to $g(y)>g(x).$ The forgotten portion of costs to
move, $\left( 1-\xi (x)\right) C(x,y)$ is traded against an expected
increase of the speed of improvement, because less exploration around is
required to estimate and satisfy the worthwhile-to-move constraint. The
ratio $1-\xi (x)\geq 0$ is the chosen sacrificing rate, while its complement 
$\xi (x)$ is the non sacrificing rate at $x\in X.$

Define the \textquotedblleft worthwhile-to-move'' relationship without
opportunity costs
\begin{equation}  \label{Punc1}
x\in X\longmapsto W(x)=\left\{ y\in X: \mbox{  } A(x,y)\geq
\xi(x)C(x,y)\right\}
\end{equation}
\noindent and the explored portion of it $W^{E}(x)=W(x)\cap E\left( x,r(x)\right).$ 
 This relationship is clarified by considering both the respective
lengths of the periods $t(x,y)$ and $t(y)$, and the advantages and costs to
move:
\begin{equation}  \label{Punc2}
A(x,y) \geq \xi (x)C(x,y)
\end{equation}
\begin{equation}  \label{Punc3}
\Longleftrightarrow \mbox{  } t(y)\delta (x)\left( g(y)-g(x)\right) \geq \xi
(x)t(x,y)c(x,y)=\xi(x)e(x,y)d(x,y).
\end{equation}
\begin{equation}  \label{Punc4}
\Longleftrightarrow \mbox{  }g(y)-g(x)\geq \theta (x,y)d(x,y)
\end{equation}
\noindent where
\begin{equation}  \label{Punc5}
\theta (x,y)= \frac{\xi (x)e(x,y)}{t(y)\delta (x)} >0,
\end{equation}
\noindent if $t(y)>0$ and $\delta (x)>0.$ The worthwhile-to-move
relationship becomes
\begin{equation}  \label{Punc6}
W(x) = \left\{ y\in X: \mbox {  } g(y)-g(x)\geq \theta (x,y)d(x,y)\right\}
\subset X.
\end{equation}
A worthwhile-to-move step is such that the ratio between the advantages to
move and the distance to move, $\left( g(y)-g(x)\right) /d(x,y)$ is higher
than the acceptable transition rate $\theta (x,y)>0$.

Consider the case where opportunity costs $O(x,y)=t(x,y)g(x)$ are included
in the comparison between advantages and costs to move. The
worthwhile-to-move relationship becomes $A(x,y)\geq \xi (x)C(x,y) + \xi
(x)O(x,y).$ Introducing the speed of moving $v(x,y)=d(x,y)/t(x,y)$ gives the
previous situation. In this case, the time spent to move is $%
t(x,y)=v(x,y)d(x,y)$ and opportunity costs to move are $%
O(x,y)=g(x)d(x,y)/v(x,y).$ The acceptable transition ratio becomes $\theta
(x,y)=\xi (x)\left[ e(x,y)+(g(x)/v(x,y))\right] /\left[ t(y)\delta (x)\right]%
.$

\subsection{ Topology and Tychastics.}

The main difficulty of decision-making  in a complex environment
is to manage the transition toward a final goal, balancing
between intermediate state-dependent advantages and costs to move. During
the transition, we acknowledge our ignorance of what the agent is exactly
doing, when costs to move are important to consider, a particular case of
the \textquotedblleft satisficing-by-rejection principle'' of Martinez-Legaz
et al. (2002), and Soubeyran (2006).

The agent rejects the moves that require too many intermediate sacrifices,
 not \textquotedblleft worthwhile - to-move'', because advantages to
move are lower than a certain proportion of costs to move. The agent rejects
alternatives which are perceived as too much sacrificing. Inequalities or
set inclusions are convenient to uncertainty with no knowledge of
probability distribution (Aubin, 2005).

At each step, the first problem is ignorance. Exploration concerns knowledge
acquisition around a given state (structured information),
instead of information acquisition. Exploration helps the agent to discover
his taste locally, with reference to the current state, and to better know
his feelings about advantages and costs to move, to estimate his resilience
to effort, to build his current ambition. It raises the question of how to
choose the radius of the exploration set $r(x)\geq 0,$ and of how much to
spend for exploring around the current state. Is it worth exploring to know
if it is whorthwhile to move? When, after having explored around a given
state, is it worthwhile to move? To decide \textit{ex ante} how much to
explore starts an infinite regression which is difficult to cut. Heuristics
of exploration can help do that approximatively. To escape from this
regression problem, we assume a constant size of the exploration set, $%
r(x)=r>0$. If the agent follows a \textquotedblleft worthwhile-to-move''
process, he will finally enter into a clairevoyance ball and optimize. We
also assume that the agent explores more or less, depending on his
motivation to change.

\smallskip

Based on a \textquotedblleft worthwhile-to-move'' behavior as the reference,
we suggest a classification of behaviors according to the following four
categories (they all are special cases of the last category):

a) \textit{Low goal-oriented behaviors} including

\begin{itemize}
\item Improving: $g(y)>g(x)$;
\item \textquotedblleft Worthwhile-to-move'', \textquotedblleft not too many
sacrificing'' behaviors:
$A(x,y)\geq \xi (x)C(x,y);$
\end{itemize}

b) \textit{High goal-oriented behaviors} including
\begin{itemize}
\item \textquotedblleft Improving enough'': $g(y)\geq g(x)+\varepsilon (x),%
\mbox{ } \varepsilon (x)>0$ which is equivalent to \textquotedblleft
intermediate satisficing'': $g(y)\geq \widetilde{g}(x)=g(x)+\varepsilon (x)$;
\item \textquotedblleft Intermediate satisficing with not too much
sacrificing'':
$A(x,y)\geq \xi (x)C(x,y)$ and $g(y)\geq \widetilde{g}(x).$
\end{itemize}

The \textquotedblleft not too many intermediate sacrificing'' constraint
defined in Eq.(\ref{WTMset1}) reflects the presence of costs to move.

In all cases, because of some ignorance, there is an adjoint state-dependent
exploration process $E(.):$ $x\in X\longmapsto E(x,r(x))$ devised to
discover the state-dependent advantages and costs to move. At each step,
before moving, the agent explores a portion of the worthwhile-to-move set, $%
W^{E}(x).$

\smallskip

The model is made of three blocks:

i) \textit{Goal setting block}: $y\in I\left[ x,\varepsilon (x)\right]
=\left\{ y\in X: \mbox{ } g(y)\geq g(x)+\varepsilon (x)\right\};$

ii) \textit{Moving block} during the transition, $y\in W(x);$

iii) \textit{Exploration block}: $y\in E\left( x,r(x)\right) .$

At each static period we introduce intermediate payoffs. The choice of the
length of each exploitation phase partially determines the size of each
intermediate payoff as an endogenous intermediate goal.
During static periods the agent exploits the benefits or the net utility $%
g(x)$ coming from the repetition of the present alternative $x$, then
explores around for an alternative $y\in E(x,r(x))$\ which improves enough.
During dynamic periods the agent moves only.

\section{Incrementalism through Enclosing}

\subsection{Using an Enclosing Heuristic to Enclose the \textquotedblleft
Worthwhile to Move'' Inclusion}

We show how the agent manages the \textquotedblleft worthwhile-to-move''
punctuated transition, using an enclosing heuristic. At each step, the agent
sets upper and lower bounds on several control variables. He encloses these
variables within intervals.

The maximum duration of exploitation is $0<t(y)\leq \overline{t}$, the
minimum per unit of time effort of moving $e(x,y)\geq \underline{e}>0$ and
the minimum rate of non sacrificing $\xi (x)\geq \underline{\xi }>0.$
Maximum weights over satisfaction and deception are $\delta (x)\leq 
\overline{\delta }.$ The acceptable transition ratio $\theta (x,y)>0$
(defined in Eq.(\ref{Punc5}) section 6) is higher than a satisficing
transition ratio $\theta =(\underline{\xi }\underline{e})/\left( \overline{t}%
\overline{\delta }\right) >0.$

In the presence of opportunity costs to move, if the speed of moving is not
too high, $0\leq v(x,y)\leq \overline{v}$, the satisficing transition ratio
is $\theta (x)=\underline{\xi }(\underline{e}+g(x)/\overline{v})/\left( 
\overline{t}\overline{\delta }\right) \geq \theta =\underline{\xi }(%
\underline{e}+g(x_{0})/\overline{v})/\left( \overline{t}\overline{\delta }%
\right) $ because along an improving trajectory $g(x)\geq g(x_{0}).$ In both
cases, following a worthwhile-to-move process, for any consecutive steps $x$
and $y$, the following inequality 
\begin{equation*}
g(y)-g(x)\geq \theta d(x,y)
\end{equation*}
holds.
Let $S(x)=\left\{ y\in X: \mbox{ } g(y)-g(x)\geq \theta d(x,y)\right\}
\supseteq W(x)$ be the enclosing worthwhile-to-move set. The enclosing
heuristic is powerful because it defines a pre-order on the state space $X$.
The enclosing worthwhile-to-move relationship $x\in X\longmapsto S(x)\subset
X$ is reflexive, $x\in $ $S(x)$ for all $x\in X$ and transitive.

If the distance effort of moving $e(x,y)$ is a strictly increasing function
of the speed of moving, $e(x,y) =\eta \left( v(x,y)\right)$, and if the
speed of moving is higher than a strictly positive minimum level, $%
v(x,y)\geq \underline{v}>0,$ the effort per unit of distance $e(x,y)$ will
be higher than a strictly positive level $e(x,y)\geq \underline{e} = \eta (%
\underline{v}) >0 $ for all $x,y\in X.$

\subsection{More on Incrementalism: Convergence in Finite Time}

Each period consists of a static period $\left[ x\right] $ of
exploration-exploitation of duration $h(x),$ where $x=x_{n}$, and of a
dynamic moving period $\left[ x,y\right] $ from $x$ to a new improving state 
$y=x_{n+1}$ of duration $t(x,y)$. A new period of exploration-exploitation $%
\left[ y\right] $ starts. At each exploration-exploitation period, the time
spent for exploitation is $t(x)=\alpha (x)h(x)$, and the time spent for
exploration is $\left( 1-\alpha (x)\right) h(x),$ with $0\leq \alpha (x)\leq
1.$

The total time spent for exploration-exploitation and moving along a
worthwhile-to-move and satisficing trajectory is
\begin{equation}  \label{CFT1}
T=\sum_{n=0}^{+\infty }\left( h(x_{n})+t(x_{n},x_{n+1})\right).
\end{equation}
If the speed of moving is lower bounded $v(x,y)\geq v>0$, the relationship
between speed and distance implies that $d(x,y)=t(x,y)v(x,y)\geq t(x,y)v$.
Subsequently,
\begin{equation}  \label{CFT2}
\sum_{n=0}^{+\infty }t(x_{n},x_{n+1})\leq \frac{1}{v}
\sum_{n=0}^{+\infty }d(x_{n},x_{n+1})<+\infty ,
\end{equation}
\noindent this last inequality being a consequence of the \textquotedblleft
worthwhile-to-move'' principle (see subsection 2.2).

If the agent exploits a fraction of time $\alpha (x)\geq \underline{\alpha }%
>0$ greater than or equal to a given minimum level $\underline{\alpha }$ $>0$
and if the total time spent for exploitation is finite, $\left(\sum_{n=0}^{+%
\infty }t(x_{n})<+\infty \right)$ then:
\begin{equation}  \label{CFT3}
\sum_{n=0}^{+\infty }h(x_{n})=\sum_{n=0}^{+\infty }\frac{t(x_{n})}{%
\alpha(x_{n})} < \frac{1}{\underline{\alpha}} \sum_{n=0}^{+\infty
}t(x_{n})<+\infty .
\end{equation}
Adding the two inequalities (\ref{CFT2}) and (\ref{CFT3}) and using Eq.(\ref%
{CFT1}) gives the convergence in finite time.

\section{Instrumentalism: Shrinking of 
Worthwhile-to-Move Behaviors which \textquotedblleft Improve Enough''}

Consider a more goal-oriented \textquotedblleft worthwhile-to-move'' process
where, in an inertia context, the agent wants both, to follow a
\textquotedblleft worthwhile-to-move'' transition 
\begin{equation*}
A(x,y)\geq \xi (x)C(x,y)
\end{equation*}
where the agent must compensate intermediate costs to move by intermediate
advantages to move, in such a way that advantages to move are higher than a
given fraction of costs to move; and to \textquotedblleft improve enough'': 
\begin{equation*}
g(y)-g(x)\geq \varepsilon (x).
\end{equation*}
We shall show that, if
\begin{itemize}
\item the instantaneous utility function is upper bounded,
\item the agent is able to enclose the \textquotedblleft
worthwhile-to-move'' inclusion $y\in W(x)$ in an enclosing inclusion $y\in
S(x),$
\end{itemize}
\noindent then, such a \textit{\textquotedblleft worthwhile-to-move''} and
\textit{\textquotedblleft temporary satisficing''} process shrinks. The radius 
 $\rho(S(x_{n})) = \sup_{y\in S(x_n)} d(x_n, y)$ of the enclosing inclusion converges to zero. Hence the
enclosing inclusion and the \textquotedblleft worthwhile-to-move'' inclusion
shrink.

\subsection{\textquotedblleft Improving Enough'' Worthwhile-to-Move Behaviors%
}

Consider first a goal-oriented setting process which is not
\textquotedblleft worthwhile-to-move''. This defines an intermediate
satisficing process where the agent wants to \textquotedblleft improve
enough'' (Soubeyran, 2006). At each step, the agent at $x$ sets a new
aspiration level $\widehat{g}(x)>g(x)$, then sets an adjoint satisficing
level $\widetilde{g}(x), \mbox{ } g(x)<\widetilde{g}(x)<\widehat{g}(x),$ and
tries to reach it. He must explore in an exploration set $E\left(
x,r(x)\right) \subset X$ around $x$ to find some $y$ such that $g(y)\geq 
\widetilde{g}(x)$. The agent adapts his aspiration level. Let $\overline{g}%
=\sup \left\{ g(y): \mbox{ } y\in X\right\}$ be the finite supremum of the
upper bounded utility function $g(.).$ The agent sets a feasible aspiration
level, $\widehat{g}(x)\leq \overline{g}<+\infty .$ Otherwise, he does not
reach it and, sooner or later, will have either to further explore around
the current state or to relax his aspiration level.

If the agent follows both a \textquotedblleft worthwhile-to-move'' $y\in
W(x) $ and an intermediate satisficing process and if he has succeeded to
enclose it in the enclosing process $y\in S(x)$, by exploration around the
current state, the agent can discover locally his utility function as well
as the enclosing inclusion
\begin{center}
$S(x)=\left\{ y\in X: \mbox{ } g(y)-g(x)\geq \theta d(x,y)\right\}.$
\end{center}
\noindent Let
\begin{equation}  \label{WTM-IIEB1}
s(x)=\sup \left\{ g(y): \mbox{ } y\in S(x)\right\} \leq \overline{g}<+\infty
\end{equation}
\noindent be the highest unknown aspiration level that the agent can reach
in the enclosing inclusion. Let
\begin{equation}  \label{WTM-IIEB2}
\widehat{g}(x)=g(x)+p(x)\left( s(x)-g(x) \right) \mbox{; with } p(x)>0
\end{equation}
\noindent be the unknown relationship between the aspiration level $\widehat{%
g}(x)$ and $s(x)$. If starting from $x$, and exploring around, the agent can
find a \textquotedblleft worthwhile-to-move'' state $\ y\in S(x)\subset W(x) 
$ which improves enough, or
\begin{equation}  \label{WTM-IIEB3}
g(y)-g(x)\geq q(x)\left( \widehat{g}(x)-g(x)\right) =\varepsilon (x)
\end{equation}
\noindent where the rate of need reduction $q(x)\in \left] 0,1\right[ $ is
greater than or equal to a minimum level, $0<\underline{q}\leq q(x)<1$ for
all $x\in X$, then, from $\widehat{g}(x)-g(x)=p(x)\left( s(x)-g(x)\right),$
the intermediate satisficing level is
\begin{eqnarray}  \label{WTM-ISL}
\widetilde{g}(x)&=& g(x)+q(x)\left( \widehat{g}(x)-g(x)\right) \\
&=& g(x) + p(x)q(x)\left( s(x)-g(x)\right) \\
&=& g(x) + \sigma (x)\left( s(x)-g(x)\right)
\end{eqnarray}
\noindent where $0<$ $\underline{\sigma }<\sigma (x)=p(x)q(x)<1.$ The
\textquotedblleft improving enough'' condition is defined by
\begin{equation}  \label{WTM-IIEB4}
g(y)-g(x)\geq p(x)q(x)\left( s(x)-g(x)\right) =\sigma (x)\left(
s(x)-g(x)\right).
\end{equation}
Take $0<$ $\underline{\sigma }<\sigma (x)=p(x)q(x)<1$, which is always
possible, although the agent must explore enough to manage to do that. This
requires to set a high enough aspiration level $\widehat{g}(x)$ with respect
to its highest feasible level $s(x)$ (the rate $p(x)$ must be high enough, $%
0<\underline{p}\leq p(x)$), and to try to fill a \textquotedblleft large
enough'' fraction $1>q(x)\geq \underline{q}$ $>0 $ of the aspiration gap $%
\widehat{g}(x)-g(x)>0$, such that $0<$ $\underline{\sigma }<\sigma
(x)=p(x)q(x)<1,$ where $0<$ $\underline{\sigma }=\underline{p}\underline{q}%
<1.$ This means that if the agent chooses a small rate $q(x)$, he can choose
a large rate $p(x)$. This is the case if he has chosen a high enough
aspiration level $\widehat{g}(x).$

Starting from any $x$, an \textquotedblleft improving enough'' and
\textquotedblleft worthwhile-to-move'' state $y\in S(x)\subset W(x),$ such
that Eq.(\ref{WTM-IIEB4}) is satisfied with $0<\underline{\sigma }<\sigma
(x)<1$, is always possible because $s(x)$ is a supremum.

\subsection{Shrinking}

The inequality  $s(x)-g(x)\geq g(y)-g(x)\geq \theta d(x,y)$ for all 
$y\in S(x)$ implies that
\begin{equation}  \label{Shrink1}
s(x)-g(x)\geq \theta \mbox{ }\rho(S(x)).
\end{equation}
\noindent We now show how the worthwhile-to-move process shrinks.
We assume that the agent can find some $x_{n+1}\in S(x_{n})\subset W(x_{n})$
which improves enough: 
\begin{equation}  \label{Shrink2}
g(x_{n+1})-g(x_{n})\geq \underline{\sigma }\left(s(x_{n})-g(x_{n})\right).
\end{equation}

\medskip

\noindent \textbf{The Shrinking Proposition:} \textit{If \newline
\indent i) the utility function is upper bounded,\newline
\indent  ii) the agent has succeeded in defining and enclosing his worthwhile-to-move
process,\newline
\indent  iii) the agent explores enough around to be able to find $x_{n+1}\in
S(x_{n})\subset W(x_{n})$ which improves enough, or such that $%
g(x_{n+1})-g(x_{n})\geq \underline{\sigma }\left[ s(x_{n})-g(x_{n})\right]
\geq 0$\\
\noindent then, the worthwhile-to-move process shrinks, which means that the
radius $\rho(W(x_{n}))$ of the worthwhile-to-move set converges to zero.}

\smallskip

\textit{Proof:} By assumption iii) and using Eq.(\ref{Shrink1}) 
\begin{eqnarray}
g(x_{n+1})-g(x_{n}) & \geq & \underline{\sigma }\left[s(x_{n})-g(x_{n})%
\right] \\
& \geq & \underline{\sigma } \theta \mbox{ }\rho(S(x_{n})).
\end{eqnarray}

We saw that a worthwhile-to-move process is such that $g(x_{n})%
\longrightarrow g^{\ast }.$ This implies that $\rho(S(x_{n}))\longrightarrow 0$
as $n\rightarrow +\infty. \Box$

\smallskip Instrumentalism is a key feature of this process: the goal
changes along the process, because the unsatisfied needs $\widehat{g}%
(x_{n})-g(x_{n})$ decrease each time.

\subsection{Convergence to a Rest point: Stopping with no Residual
Frustration}

When does the agent prefer to stop moving? In such cases, his aspiration gap
must vanish, or $\widehat{g}(x^{\ast })-g(x^{\ast })=0,$ with no residual
frustration.

\smallskip

Behavioral Rest Points: As seen before, a\textbf{\ }performance $x^{\ast
}\in X$ is said to be a behavioral rest point if, for any $y\in X$ with $%
y\neq x^{\ast }$, it is not worthwhile to move from $x^{\ast }$ to $y$. This
is equivalent to say that $x^{\ast }$ is a rest point element of the
\textquotedblleft worthwhile-to-move'' relationship $x\in X$ $\longmapsto
W(x)\subset X,$ or $W(x^{\ast })=\left\{ x^{\ast }\right\}.$ The agent has
no further incentive to move. Thus, $x^{\ast }$ is a behavioral rest point
iff for any $y\neq x^{\ast },$ $A(x^{\ast },y)<\xi (x^{\ast })C(x^{\ast },y)$%
, or with opportunity costs $A(x^{\ast },y)<\xi (x^{\ast })\left( C(x^{\ast
},y)+O(x^{\ast },y)\right) .$ This is the case if $S(x^{\ast })=\left\{
x^{\ast }\right\} $ because $S(x)\supseteq W(x)$ and $x\in W(x)$ for all $%
x\in X.$

In terms of instantaneous advantages and costs to move, $x^{\ast }\in X$ is a behavioral
rest point if, for any $y\neq x^{\ast },$
\begin{center}
$g(y)-g(x^{\ast })<\theta (x^{\ast },y)d(x^{\ast },y).$
\end{center}
In previous sections 2.2 and 8.2, we showed that a \textquotedblleft
worthwhile-to-move'' process which \textquotedblleft improves enough,''
namely $\left\{x_{n+1}\in W(x_{n}), \mbox{  }n\in{\mathbb{N}} \right\}$
converges toward some limit $x^{\ast }\in X$, and that the radius $%
\rho(S(x_{n}))$ of the enclosing inclusion $S(x_{n})\supseteq W(x_{n}), \ n\in{%
\mathbb{N}} ,$ goes to zero. But this does not imply that $\rho(S(x^{\ast })=0,$
which is equivalent to $S(x^{\ast}) = \left\{ x^{\ast}\right\}.$ When is it
the case? When the agent stays at $x^{\ast }$ rather than moves again? This will
support the existence of behavioral rest points. The upper semicontinuity of
the utility function $g(.)$ is a sufficient condition. The agent follows a
punctuated dynamic of temporary routines to reach a permanent routine. For
the sake of clarity, we give a reduced form of our model.

\paragraph{\textbf{A Reduced Form of the \textquotedblleft
Worthwhile-to-Move'' Model. }}

A \textquotedblleft worthwhile-to-move'' and \textquotedblleft improving
enough'' process has the following reduced form:

i) Starting from $x=x_{n}\in X$ and $g(x)\in\mathbb{R }$, set an
intermediate aspiration level $\widehat{g}(x)>g(x) $ and set a feasible
intermediate satisficing level $\widetilde{g}(x) $ such that $g(x)<%
\widetilde{g}(x)<\widehat{g}(x).$

ii) Around the current state $x\in X$, define a subset of local acceptable
transitions $y\in W(x),$ which determines a worthwhile-to-move behavior 
$W(.):$ $x\in X\longmapsto W(x)\subset X,$ $x\in W(x)$ at each step.

Problem:

iii) Explore around the current state $x\rightarrow y\in X,$ in the
exploration set $E\left(( x,r(x)\right) $, choosing the size $r(x)\geq 0$ of
the exploration set.

iv) Find $y=x_{n+1}\in W(x)\cap E\left( x,r(x)\right) =W^{E}(x)\subset X$
such that $g(x)<\widetilde{g}(x)\leq g(y)\leq \widehat{g}(x).$

v) If $g(y)<\widehat{g}(x)$, start from $y\in W^{E}(x)$ and set a new
intermediate goal $\widehat{g}(y)\leq \widehat{g}(x);$

and so on.

vii) Stop at $x^{\ast }\in X,$ when $W(x^{\ast })=\left\{ x^{\ast }\right\},$%
 \ when it is no longer worthwhile to move: the satisficing paradox is
solved. In that case, the agent has no residual frustration feeling, because
his aspiration gap is zero: $\widehat{g}(x^{\ast })-g(x^{\ast })=0.$

\medskip

\noindent \textbf{The Worthwhile-to-Move Theorem.}
\textit{Assume that\\
\indent i) the state space $X$ is a metric space with metric $d ,$\\
\indent ii) the instantaneous utility function $g(.)$ is upper bounded,\\
\indent iii) the agent uses a heuristic enclosing the worthwhile-to-move inclusion $y\in W(x)$
\begin{center}
 $W(x)=\left\{ y\in X: \mbox{  } A(x,y)\geq \xi (x)C(x,y)\right\} $
\end{center}
 in the inclusion  $y\in S(x)$,
\begin{center}
$S(x) = \left\{ y\in X: \mbox{  } g(y)-g(x)\geq \theta d(x,y)\right\}, \mbox{  } \theta >0.$ 
\end{center}
This means that the agent limits the control variables: he sets
a maximum duration of exploitation $0<t(y)\leq \overline{t}$, a minimum
effort of moving $e(x,y)\geq \underline{e}>0$ and a minimum rate of non
sacrificing $\xi (x)\geq \underline{\xi }>0.$ He also puts maximum weights
on satisfaction and deception, and $0\leq \delta (x)\leq \overline{\delta }.$
Then, the acceptable transition ratio is greater than or equal to a strictly
positive level, $\theta (x,y)\geq \theta >0,$ where the minimum acceptable
transition ratio is $\theta =(\underline{\xi }\underline{e})/\left( 
\overline{t} \overline{\delta }\right) >0.$\\
Then,\\
a) The worthwhile-to-move, but, for the moment, not satisficing process $%
y\in W(x)$ can be enclosed in the nested enclosing process $y\in S(x)$:
\begin{center}
$W(x)\subset S(x)$ for all $x\in X.$
\end{center}
b) Let $x_{n+1}\in W(x_{n}) $ and $x_{n+1}\in S(x_{n}), \mbox{ 
} n\in{\mathbb{N}} , \mbox{  }x_{0}\in X$ given, be the worthwhile-to-move
and the enclosing process. If the state space is complete, the
worthwhile-to-move process converges, $x_{n}\longrightarrow x^{\ast }\in X$
as $n\rightarrow +\infty $ whatever the starting state $x_{0}\in X$. It
converges in finite time if, along the process, the total time spent for
exploitation is finite and the speed of moving $v(x,y)=d(x,y)/t(x,y)$ is
greater than or equal to a strictly positive level $v>0.$\\
c) Moreover, if the worthwhile-to-move process improves enough, which is the
case if the agent \textquotedblleft explores enough'' around, at each step,
then, the process not only converges, but also shrinks: $\rho(S(x_{n}))%
\longrightarrow 0.$\\
If the per unit of time utility function is upper semicontinuous, then, the
limit state $x^{\ast }\in X$ is a behavioral rest point: $S(x^{\ast
})=\left\{ x^{\ast }\right\} \Longrightarrow $ $W(x^{\ast })=\left\{ x^{\ast
}\right\}.$\newline
The agent stops at $x^{\ast }$ with no residual frustration: $\widehat{g}%
(x^{\ast })=g(x^{\ast }).$}

\smallskip

\textit{Proof}: Only the last point c) needs proving.

Each enclosing set $S(x)\subset X$ is closed, this is a consequence of the
upper semicontinuity of $g$. The enclosing inclusion is nested. If $%
x_{n+1}\in S(x_{n})$ for all $n\in{\mathbb{N}},$ and $p\geq n $ then, $%
S(x_{0})\supset S(x_{1})\supset ...\supset S(x_{n})\supset ...\supset
S(x_{p})\supset ... $ Thus, $x_{p}\in S(x_{n})$ for all $p\geq n, $ for any
given $n\in{\mathbb{N}}$. From $x_{p}\longrightarrow x^{\ast } \mbox{ as }
p\longrightarrow +\infty $, $x^{\ast }\in S(x_{n})$ for all $n\in{\mathbb{N}}%
,$ because $S(x)$ is closed for all $x\in X.$ By transitivity of $S$, $%
S(x^{\ast })\subset S(x_{n})$ for all $n\in \mathbb{N}$, which  implies, by the
Shrinking Proposition, that $0\leq \rho(S(x^{\ast}))\leq
\rho(S(x_{n}))\longrightarrow 0.$ This gives $\rho(S(x^{\ast }))=0\Longrightarrow
S(x^{\ast })=\left\{ x^{\ast }\right\} $ because $x\in S(x)$ for all $x\in
X. $ Thus, $W(x^{\ast })=\left\{ x^{\ast }\right\} $ because of the
enclosing heuristic $W(x)\subset S(x)$ for all $x\in X.$

\subsection{Link with Ekeland's $\protect\varepsilon $-Variational Principle}

As a striking application of the \textquotedblleft worthwhile-to-move''
theorem, one obtains a cognitive version of  Ekeland's variational
theorem (Aubin and Ekeland, 1984; Attouch and Soubeyran, 2006). This theorem
(Ekeland, 1974) was originally a regularization theorem devoted to
ill-behaved optimization problems, through approximate optimization. We give
its \textquotedblleft supremum formulation\textquotedblright\ in order to
follow the tradition in economics which is usually concerned with
maximization or approximate maximization problems. It concerns the case
where the function $g(.)$ possibly has no maximum on $X$.

\smallskip

\noindent \textbf{Ekeland's Theorem}: \textit{Let $X$ be a complete metric space with distance 
$d$. Let $g(. ):x\in X\longmapsto g(x)\in\mathbb{R }$ be an upper bounded
and upper semicontinuous function. Let $\overline{g}= \sup _{x\in X} g(x)
<+\infty $, $\mbox{ }\theta >0$, and $\varepsilon >0.$ Then, for any $%
x_{0}\in X$ such that $n(x_{0})=\overline{g}-g(x_{0})\leq \varepsilon $,
there exists a $x^{\ast }\in X$ such that i), ii) and iii) are satisfied:\\
\indent i) $g(x^{\ast })\geq g(x_{0}),$\\
\indent  ii) $\varepsilon \geq \theta d(x_{0},x^{\ast }),$\\
\indent iii) $g(y)-\theta d(x^{\ast },y)<g(x^{\ast })$ for all $y\neq x^{\ast },$
or 
$x^{\ast }\in \arg \max \left\{ g(y)-\theta d(x^{\ast },y): \mbox{ }y\in
X\right\}.$}

\smallskip

The condition $n(x_{0})=\overline{g}-g(x_{0})\leq \varepsilon $ means that
the initial need or frustration $n(x_{0})$ is less than a given $\varepsilon
>0.$ Statement i) means that the final position $x^{\ast }\in X$ is
improving with respect to the initial position $x_{0}$. Statement ii) tells
us that the highest possible advantage $\varepsilon $ to move from the
initial to the final position is greater than the cost to move $\theta
d(x_{0},x^{\ast })$ from $x_{0}$ to $x^{\ast }$. It is worthwhile to move
from $x_{0}$ to $x^{\ast }.$ The last statement iii) says that $x^{\ast }$
is a maximum of the approximate function $y\in X\mapsto g(y)-\theta
d(x^{\ast },y)\in\mathbb{R}.$ It tells us that $x^{\ast}$ is a behavioral
rest point because it is not worthwhile to move from $x^{\ast }$ to any
different position $y\neq x^{\ast }.$

\subsection{From Inertia Inefficiency to Behavioral Rest Point}

Why does a \textquotedblleft worthwhile-to-move'' transition process end in
an inefficient rest point exhibiting a large inefficiency gap $\overline{g}%
-g(x^{\ast })>0$? The answer is: because of inertia inefficiencies, defined
as resistance to change. Their modelling involves \textquotedblleft costs to
move'', including costs to build and to sustain motivation, to set goals, to
explore locally, and to move.

\smallskip

What does make an agent behave inefficiently? What does generate local
actions, small steps to move, a long transition time, large intermediate
sacrifices, large per unit of time costs to move, a long time spent in moving,
small intermediate improvements and low intermediate advantages to move, a
short intermediate exploration and exploitation? What does generate
premature convergence, local optimum, a low limit level of the final goal $%
g^{\ast}$, a large inefficiency gap $\overline{g}-g^{\ast }>0,$ convergence
in a very long time, a low speed of convergence, an irregular punctuated
dynamic, some permanent frustration? How can we explain the emergence of
routines, habits, and rest points? How do the ambivalent aspects of habits
and routines, both positive and negative, reach a balance? How can inertia
lead to inefficient behaviors?

Our answer is given through the two main inequalities which  the agent has to manage at each step (we set
$x=x_{n}$ and $y=x_{n+1}$)

\smallskip

i) \textquotedblleft not too many intermediate sacrifices''
\begin{equation}  \label{II1}
g(y)-g(x)\geq \theta (x,y)d(x,y)
\end{equation}
ii) \textquotedblleft improving enough''
\begin{equation}  \label{II2}
g(y)-g(x)\geq q(x)\left[ \widehat{g}(x)-g(x)\right] = \varepsilon (x)>0%
\mbox{, or}
\end{equation}
iii) choose to stop, setting $\varepsilon (x)=0.$

The agent manages the \textquotedblleft worthwhile-to-move'' inequality Eq. (\ref{II1}) by choosing the \textquotedblleft adjusted non sacrificing
ratio,''  using an enclosing heuristic
\begin{equation}  \label{Ansr}
\theta (x,y)=\frac {\xi (x)e(x,y)}{t(y)\delta (x)} \geq \underline{\theta }>0
\end{equation}
and manages the \textquotedblleft improving enough'' inequality Eq. (\ref{II2})
 by choosing the satisficing level $\mbox{ } \widetilde{g}(x)$ and the
improvement rate $q(x). $

In this context, behavioral inefficiencies come from high goal-setting,
exploration costs, and moving costs. Inertia generates a lot of
inefficiencies if the enclosing process gives a high lower bound $\underline{%
\theta }>0$, and if $\widetilde{g}(x)$ and $q(x) $ are low. Behavioral
inefficiencies and premature convergence come from limited needs and strong
inertia with not enough intermediate sacrificing, coming from too high a
\textquotedblleft non sacrificing rate'' $\xi (x)$, too short a horizon $h(y), $ too strong a preference for the present, too low an exploitation time 
$t(y)=\alpha (y)h(y)$, too high an effort to move $e(v(x))$ (coming from too
high a speed of moving), too low a psychological factor $\delta (x)$, too
high local exploration costs $K(x)$.

Then, how, at each step, are determined the \textquotedblleft adjusted non
sacrificing ratio'' $\theta (x,y)$, the intermediate satisficing level $%
\widetilde{g} (x)$, the rate of \textquotedblleft improving enough'' $q(x),$
and the exploration expenditures $K(x)$?

\smallskip

Assume that the agent has limited resources and a finite supremum $\overline{%
g}<+\infty $ which is the highest utility he can reach. Starting from the
state $x\in X,$ let $\widehat{g}(x)$ be the aspiration level of the agent,
which is an estimation of his unknown supremum utility. The lower his
incentive to change, the lower his initial unsatisfied needs $\widehat{g}%
(x)-g(x)>0.$ The agent can over- or under-estimate the supremum $\overline{g}$. 
This depends on self-esteem and degree of optimism. If his aspiration
level is higher than his feasible ambition, $\widehat{g}(x)>\overline{g},$
the agent will sooner or later, after trials and errors, have to relax his
aspiration level, so as to make it feasible, say $\widehat{g}(x)\leq 
\overline{g}$. To fulfill this goal-setting task, the agent can spend time
and money to imitate successful agents or ask a coach for advice. The agent
can try to reach his intermediate aspiration level in several steps. He can
set an intermediate satisficing level, an intermediate goal $\widetilde{g}%
(x) $ between his present utility level and his intermediate aspiration
level such that $g(x)<\widetilde{g}(x)<\widehat{g}(x).$ Then he can try to
find an intermediate satisficing state $y\in W(x)$ such that $g(x)<%
\widetilde{g}(x)\leq g(y)\leq \widehat{g}(x).$

Consider the intermediate goal-setting inequality 
\begin{equation*}
g(y)-g(x)\geq q(x)\left( \widehat{g}(x)-g(x)\right) =\varepsilon (x)>0.
\end{equation*}
To set such an \textquotedblleft improving enough'' inequality requires to
be able to set an intermediate aspiration level $\widehat{g}(x)=g(x)+p(x)%
\left(s(x)-g(x)\right) $ and a related intermediate satisficing level $%
\widetilde{g}(x)=g(x)+\varepsilon (x)=g(x)+q(x)\left(\widehat{g}(x)-g(x)%
\right) .$ The equation $\widehat{g}(x)-g(x)=p(x)\left( s(x)-g(x)\right)$
 defines the intermediate satisficing level
\begin{equation}
\widetilde{g}(x)=g(x)+p(x)q(x)\left( s(x)-g(x)\right) = g(x)+\sigma
(x)\left( s(x)-g(x)\right),
\end{equation}
\noindent where $\mbox{ } 0<$ $\underline{\sigma }<\sigma (x)=p(x)q(x)<1.$

At each step, to set a feasible intermediate satisficing level $\widetilde{g}%
(x)\leq \overline{g}$ may be very costly. If the agent sets an unfeasible
level $\widetilde{g}(x)>\overline{g}$, he may take a long time to explore
and discover that there is no intermediate satisficing state $y$ such that $%
g(y)\geq \widetilde{g}(x).$ The agent must either further explore around the
current state $x\in X$ or must relax the intermediate satisficing level $%
\widetilde{g}(x)$, relaxing or not his intermediate aspiration level $%
\widehat{g}(x).$

If the agent succeeds in reaching the intermediate satisficing level, $%
g(y)\geq \widetilde{g}(x)$, he starts in goal setting process again, sets a
new intermediate aspiration level $g(y)<\widetilde{g}(y)<\widetilde{g}(x)<%
\widehat{g}(y)<\widehat{g}(x)$, relaxes this level or not. Costs of setting
intermediate goals are the time lost for trials and errors to set feasible
intermediate aspiration levels. Failure to reach an intermediate satisficing
level can lower his self-esteem and his motivation to improve. There are
also costs to set too low aspiration levels which do not lead the agent to
explore enough.

The agent may be over-satiated, when his intermediate residual frustration
feelings $\mu (x)$ are too low, or his intermediate satisfaction feelings $%
\lambda (x)$ too high, generating too low motivation for change. These terms
act positively and negatively on the intermediate goal level $\widehat{g}%
\left( x,\lambda (x),\mu (x)\right) . $ The agent will have an incentive to
\textquotedblleft improve enough,''  to set a high enough intermediate
satisficing level $\widetilde{g}(x),$ which will give him the motivation to
explore sufficiently, if his satisfaction to improve $\lambda (x)$ remains
high, or his discomfort feelings $\mu (x)$ to have unsatisfied needs remain
high during the process. Discount factors can also play a role to estimate
advantages to change and costs to change in a different way.

\smallskip

Consider the \textquotedblleft worthwhile-to-move'' inequality (the
\textquotedblleft not too many intermediate sacrifices'' condition)
\begin{equation}
g(y)-g(x)\geq \theta (x,y)d(x,y).
\end{equation}
\noindent The larger the inertia index
\begin{equation}
\theta (x,y)= \frac{\xi (x)e(x)}{(\lambda (x)+\mu (x)+\nu (x))t(y)\delta(x)},
\end{equation}
\noindent with $t(y)=\alpha (y)h(y),$ the higher the final inefficiency gap $%
\overline{g}-g^{\ast }\geq 0.$ Small intermediate steps $d(x,y)$ may be
necessary to divide very convex costs $C(x,y)$ to change. But long
intermediate steps may be necessary to improve each step enough. Too high a
speed to move $v(x,y)=d(x,y)/t(x,y)$ may be very costly. Moving may take too
long an intermediate time $t(x,y)$ to start, to move, and to stop. Enclosing
more or less can lead to more or less premature convergence. This may have
contradictory effects. A lower index of inertia $\theta (x,y)\geq \underline{%
\theta }(x)>0$ gives a larger cone where the agent can improve enough, with
not too much intermediate sacrificing, hence a more global \ process, and a
higher final utility.

Sacrificing not enough at each period ($\xi (x)$ high) may generate
premature convergence, or does not sustain intermediate motivation to reach
the final goal, because the agent does not compensate enough intermediate
costs to move by intermediate advantages to move. Other drawbacks are to
take too long a time to change, to bear too high a cost to change, or to
choose too long a time to exploit and to reach the final goal in a reasonable
time, adopting a conservative behavior, not learning enough. Furthermore, to
exploit too short a time $t(y)$ at each period discourages
motivation.

\smallskip

Consider the exploration expenditure function 
\begin{equation*}
K(x)=\tau (x)k(x) = \left( 1-\alpha (x)\right) h(x)k(x)
\end{equation*}
\noindent which represents the amount of time, energy, and money which the
agent chooses to spend to explore around a current state $x\in X$. It is
usually an increasing function of the unsatisfied needs $\widehat{g}%
(x)-g(x), $ the rate $\lambda (x)$ of contentment feelings rate and the rate 
$\mu (x)$ of discomfort to have missed the intermediate aspiration level$.$

Reference-dependent payoffs forbid to achieve exploration in a single step.
They require to explore further. Even if, starting from a given state $x$,
the agent has explored a large part of the whole space to know not only his
utility $g(y)$ but also his costs and time spent to move $C(x,y)$ and $%
t(x,y) $, he will have to explore further, starting from his new state $y\in
X$, to discover his new cost and time spent to move from $y,$ $C(y,z),$ and $%
t(y,z). $ To know before exploration how much to explore forbids
optimization. Too large or too complex a state space requires to explore
each step. The smaller the state space, the more global should the
exploration process be.

The agent may  explore not enough around the current state or too
much. He may be too shy or too uninterested to explore. \textquotedblleft
Improving enough'' requires to \textquotedblleft explore enough''. The agent
will explore enough if he has low exploration costs, or a high enough
motivation to explore, which lowers the negative feeling of exploration
expenditures, and opportunity costs, a horizon which is far away. The longer
time the agent spends exploring, the less the agent can spend in exploiting.
This may decrease his motivation to improve and the time spent to converge
may be too long.

A quick convergence may be good and bad. The speed of convergence may be
valuable, but not so much if the agent is locked in a low level
instantaneous utility. Convergence may happen earlier than for a hill
climbing process (a local search optimization with no inertia), but perhaps
too early! The process may reach too low a local optimum, or more generally
a rest point (which can be below a local maximum reached by a hill climbing
algorithm!).

\subsection{The Clairevoyance Theorem}

Consider that, at each step, the agent chooses the same radius of
exploration $r(x_{n})=r>0, \mbox{  } n\in{\mathbb{N}} .$ Our result shows
that, after some finite time, the worthwhile-to-move set lies inside the
exploration ball of constant radius, because the radius of the
worthwhile-to-move set skrinks to zero. After a finite time the agent will
optimize. This is a powerful result which shows when, under inertia and
frictions, an agent optimizes. This helps understand the degree of validity
of the \textquotedblleft as if '' hypothesis: even if agents do not
optimize, it is \textquotedblleft as if'' agents would optimize.

\section{The \textquotedblleft As if Hypothesis" Revisited: \ Local Search
and Proximal Algorithms}

In the \textquotedblleft as if hypothesis" economists make \textquotedblleft
as if" agents would optimize, even though they do not think it is the case
in many situations (Friedman, 1953).

\subsection{The \textquotedblleft As if Hypothesis" with Inertia:
\textquotedblleft A Local Search and Proximal Algorithm"}

A natural way to improve the \textquotedblleft as if hypothesis" consists of
introducing inertia costs into classical optimization programs (although
Simon does not like to add further costs to solve the bounded rationality
problem, see Gilovich and al., (2002): 559-582).

Optimization algorithms ignore inertia costs, except, in a very implicit
way, proximal algorithms (Iusem, 1995; Attouch and Teboulle, 2004; Attouch
and Bolte, 2009). We have to interpret the added regularization term which
characterizes proximal algorithms as a \textquotedblleft cost to change".
This quite simple but important interpretation completely changes our view
on the \textquotedblleft as if hypothesis". Our model shows that agents can
manage inertia \textquotedblleft as if " they use a new algorithm which is a
mixture between the two following optimization algorithms:

i) Local search algorithms: $\mbox{ } \sup \left\{ g(y): \mbox{ } y\in
E\left( x_{n},r(x_{n})\right) \right\} $

where $E\left( x_{n},r(x_{n})\right) \subset X$ is the exploration set at $%
x=x_{n}\in X,$ of size $r(x_{n})\geq 0. $ Hill climbing and simulated
annealing algorithms belong to this class of algorithms.

ii) Proximal algorithms: $\mbox{ } \sup \left\{ g(y)-\theta c(x_{n},y): %
\mbox{ } y\in X\right\}$

\noindent where $c(x,y)\geq 0$ is the regularization term which makes the
goal $g(y)$ more regular. For us, the regularization term is a cost to move,
which makes proximal algorithms satisfy the \textquotedblleft
worthwhile-to-change" relationship. In proximal algorithms, the criterion
which is to maximize can be interpreted as a net gain function, which is a
way to handle the multi-criteria problem (improve the gain function and
satisfy without too much sacrificing).

\medskip

The \textquotedblleft as if" mixture consists in solving the optimization
problem
\begin{equation}
\sup _ {y\in E\left( x_{n},r(x_{n})\right)}\left( g(y)-\theta c(x_{n},y)
\right)
\end{equation}
\noindent in order to pass from $x_{n}$ to $x_{n+1}$. This \textquotedblleft
worthwhile to change" algorithm is the \textquotedblleft Local Search and
Proximal" algorithm, LSP algorithm in short. As for classical optimization,
one can assume that, at stage $n$, the agent optimizes up to some
approximation level $\epsilon_n$,
\begin{equation}  \label{LSPA1}
x_{n+1} \in \epsilon_n-\mbox{argmax} \left\{ g(y)-\theta c(x_{n},y):\mbox{ }
y\in E\left( x_{n},r(x_{n})\right) \right\}.
\end{equation}
Starting from the current state $x_{n}\in X$ at step $n$, proximal
algorithms take the same exploration set, the whole space $E\left(
x_{n},r(x_{n})\right) =X$. From a behavioral point of view this is not
reasonable, because this means to solve a global optimization problem at
each step (the substantive case)! Hence, proximal algorithms do not consider
the recursive cost problem which is \textquotedblleft how to choose the size 
$r(x_{n})\geq 0$ of the exploration set at each step,"   \textquotedblleft how
much to explore, depending on the costs of exploration". This problem is
related to the famous \textquotedblleft effort-accuracy" trade-off related
to exploration costs and the quality of the decision (Payne and al., 1998;
Busemeyer and Diederich, 2000). To save space, we have focussed our
attention on the dynamic trade-off which balances advantages to move $%
g(y)-g(x)$ and costs to move $c(x,y).$

We have considered the simplest case of \textquotedblleft high local costs
to move'' with an exploration set of given size, where $r(x_{n})=r>0.$ The
\textquotedblleft clairevoyance theorem" shows that, in finite time, the
exploration set becomes a \textquotedblleft clairevoyance ball" where the
agent can optimize locally. In that case, the convergence of
\textquotedblleft Local Search and Proximal Algorithm" is a straight
consequence of our previous results, the agent reaches a behavioral rest
point $x^{\ast }\in X$ \ where he prefers to stay than to move, setting $%
r(x^{\ast })=0.$

We  show that the \textquotedblleft Local
Search and Proximal Algorithm" still enjoys the convergence properties in
the case of \textit{low local costs to move}.

\medskip

\noindent \textbf{Theorem: Convergence of the LSP Algorithm with Low Local Costs to
Move}

\textit{Let $X = \mathbb{R}^p$ be equipped with the Euclidean distance $d(x,y) =
\|x-y\| = (\sum _{i=1}^p (x_i - y_i)^2)^ {\frac{1}{2}}$. Let $C \subset X$
be a closed convex nonempty subset of $X$ (set of constraints, resources).%
\newline
Let $g(. ):x\in X\longmapsto g(x)\in\mathbb{R }$ be a function (gain,
utility) which satisfies i), ii) and iii):\newline
\indent i) $g$ is upper bounded on $C$, let $\overline{g} < + \infty $ be the
supremum of $g$ on $C$;\newline
\indent ii) $g$ is a smooth function (continuously differentiable);\newline
\indent iii) $g$ is quasi-concave (with convex upper level sets).\\
Given some initial data $x_0 \in X$, let $(x_n)_{n\in \mathbb{N}}$ be a
sequence defined by the Local Search and Proximal Algorithm with
clairvoyance radius $r>0$ and parameter $\theta >0$:
\begin{equation}  \label{Prox1}
x_{n+1} \in \mbox{argmax} \left\{ g(y)-\theta \|x_{n}-y\|^2: \mbox{ }y\in C,%
\mbox{ } \|y - x_{n}\| \leq r \right\}.
\end{equation}
Then, the sequence $(x_n)_{n\in \mathbb{N}}$ converges in $X$ to some $%
x_{\infty}$ which is a critical point of $g$ over $C$, namely
\begin{center}
$-\nabla g (x_{\infty}) + N_C (x_{\infty}) \ni 0$
\end{center}
\noindent where $N_C (x_{\infty})$ is the (outward) normal cone to $C$ at $%
x_{\infty}$.}

\smallskip

\textit{Proof:} By taking $y=x_n$ in Eq. (\ref{Prox1}) of $x_{n+1} $, it is
worthwhile to move from $x_n$ to $x_{n+1}$:
\begin{equation}  \label{Prox2}
g(x_{n+1}) - g(x_n) \geq \theta \|x_{n+1} - x_{n}\|^2.
\end{equation}
Summing up these inequalities and using i):
\begin{equation}  \label{Prox3}
\sum _{i=0}^{\infty} \|x_{n+1} - x_{n}\|^2 \leq \frac{1}{\theta } ( 
\overline{g} - g(x_0)) < + \infty.
\end{equation}
As a consequence
\begin{equation}  \label{Prox4}
x_{n+1} - x_{n} \longrightarrow 0 \mbox{ }\mbox{as}\mbox{ } n
\longrightarrow +\infty.
\end{equation}
Hence, for $n$ large enough, $\|x_{n+1} - x_{n}\| < r$, which implies that
the supremum in Eq. (\ref{Prox1}) is actually achieved in the interior of the
ball $\mathbb{B }(x_n,r)$ with center $x_n$ and radius $r>0$. When writing
the optimality conditions, for $n$ large enough, this exploration constraint
is not active and
\begin{equation}  \label{Prox5}
-\nabla g (x_{n+1}) + N_C (x_{n+1}) + 2\theta (x_{n+1}-x_{n}) \ni 0.
\end{equation}
The convergence of the sequence of values $(g(x_n))_{n\in \mathbb{N}}$ is
also a direct consequence of Eq. (\ref{Prox2}). It is an increasing
upper-bounded sequence. We set
\begin{center}
$g_{\infty} = \lim_{n\rightarrow +\infty} \mbox{ }g(x_n)$.
\end{center}
To prove the convergence of the sequence $(x_n)_{n\in \mathbb{N}}$ we
introduce the set
\begin{center}
$S = \{x\in C: g(x) \geq g_{\infty}\}$
\end{center}
and prove that:\newline
a) For every $a\in S, \lim_{n\rightarrow +\infty}\mbox{ } \|x_{n}-a\|^2$
exists.\newline
b) Every limit point of the sequence $(x_n)_{n\in \mathbb{N}}$ belongs to S.

Indeed, by convexity of the norm
\begin{equation}  \label{Prox6}
\|x_{n}-a\|^2 - \|x_{n+1}-a\|^2 \geq 2 \left\langle x_{n+1}-a,
x_{n}-x_{n+1}\right\rangle.
\end{equation}
By Eq. (\ref{Prox5}) there exists some $\xi_n \in N_C (x_{n+1})$ such that
\begin{equation}  \label{Prox7}
x_{n}-x_{n+1} = \frac{1}{2\theta}(-\nabla g (x_{n+1}) + \xi_n).
\end{equation}
Combining Eq. (\ref{Prox6}) and Eq. (\ref{Prox7})
\begin{equation}  \label{Prox8}
\|x_{n}-a\|^2 - \|x_{n+1}-a\|^2 \geq \frac{1}{\theta} \left\langle x_{n+1}
-a, -\nabla g (x_{n+1}) + \xi_n\right\rangle.
\end{equation}
We use the quasi-concavity assumption on $g$ in order to prove that
\begin{equation}  \label{Prox9}
\left\langle x_{n+1} - a, -\nabla g (x_{n+1}) + \xi_n\right\rangle \geq 0.
\end{equation}
To that end we consider the set
\begin{center}
$D_n = \{x\in C: g(x) \geq g( x_{n+1}) \}$
\end{center}
\noindent which is an upper level set over $C$ of the function $g$. Because
of the quasi-concavity of $g$ and of the convexity of $C$ this a closed
convex subset of $X$. By a classical geometrical argument (Rockafellar and
Wets, 2004, ch. 10)
\begin{equation}
N_ { D_n }( x_{n+1}) = - \nabla g (x_{n+1}) + N_C (x_{n+1}).
\end{equation}
Hence,
\begin{equation}
-\nabla g (x_{n+1}) + \xi_n \in N_ { D_n }( x_{n+1})
\end{equation}
\noindent and as $a \in S \subset D_n$ (recall that $g(a) \geq g (x_{n+1})$)
we obtain Eq. (\ref{Prox9}). Returning to Eq. (\ref{Prox6}) we obtain that $%
\|x_{n}-a\|^2$ is a decreasing sequence, hence converges, which proves point
a).\newline
Concerning point b), recall that $g_{\infty} = \lim_{n\rightarrow +\infty} %
\mbox{ }g(x_n)$. As $g$ is continuous, any limit point $x^{\ast }$ of the
sequence $(x_n)_{n\in \mathbb{N}}$ satisfies $g(x^{\ast }) = g_{\infty}$. As 
$C$ is a closed set and $x_n \in C$ for all $n\in \mathbb{N}$ we still have
that $x^{\ast } \in C$ at the limit. These two results imply $x^{\ast } \in
S $, which is point b).

The Opial argument consists of deducing that the whole sequence $(x_n)_{n\in 
\mathbb{N}}$ converges from a) and b). This sequence is bounded because $%
\lim_{n\rightarrow +\infty}\mbox{ } \|x_{n}-a\|^2$ exists for every $a\in S$
and $S\neq \emptyset$. If the sequence $(x_n)_{n\in \mathbb{N}}$ has two
limit points, set
\begin{equation}
x_{n_1} \longrightarrow x^{\ast }_1 \mbox{ and } x_{n_2} \longrightarrow
x^{\ast }_2 .
\end{equation}
\noindent By point b), \ $x^{\ast }_1$ and $x^{\ast }_2 $ belong to $S$. 
\newline
Using point a),\ $\lim_{n\rightarrow +\infty}\mbox{ } \|x_{n}-x^{\ast
}_1\|^2 $ \ and \ $\lim_{n\rightarrow +\infty}\mbox{ } \|x_{n}-x^{\ast
}_2\|^2$ exist. Subsequently,
\begin{equation}
\lim_{n\rightarrow +\infty}\mbox{ }\left( \|x_{n}-x^{\ast }_1\|^2 -
\|x_{n}-x^{\ast }_2\|^2\right) \mbox{  exists }
\end{equation}
\noindent which after simplification yields
\begin{center}
$\lim_{n\rightarrow +\infty}\mbox{ } \left\langle x_{n}, x^{\ast }_2 -
x^{\ast }_1 \right\rangle $ exists.
\end{center}
Specializing this result to the two subsequences $x_{n_1}$ and $x_{n_2}$
which converge respectively to $x^{\ast }_1 $ and $x^{\ast }_2 $,
\begin{equation}
\left\langle x^{\ast }_1, x^{\ast }_2 - x^{\ast }_1 \right\rangle =
\left\langle x^{\ast }_2, x^{\ast }_2 - x^{\ast }_1 \right\rangle
\end{equation}
\noindent that is, $\|x^{\ast }_1 - x^{\ast }_2\|^2 = 0$. Hence the sequence 
$(x_n)_{n\in \mathbb{N}}$ has a unique limit point and converges in $X$ to
some $x_{\infty}$.

By passing to the limit on Eq. (\ref{Prox5})
\begin{equation}
-\nabla g (x_{n+1}) + N_C (x_{n+1}) + 2\theta (x_{n+1}-x_{n}) \ni 0
\end{equation}
\noindent and by using Eq. (\ref{Prox4}) together with the smoothness of $g$
(assumption ii)) and the closedness property of the graph of the normal cone
mapping $x\mapsto N_C (x)$, we finally obtain that $x_{\infty}$ is a
critical point of $g$ over $C$, namely
\begin{equation}
-\nabla g (x_{\infty}) + N_C (x_{\infty}) \ni 0.
\end{equation}

\bigskip

In our revisited \textquotedblleft as if hypothesis", in a context of
inertia, most of our behaviors are intermediate goal-setting behaviors. It
can be seen as a mixture of a local search optimization model with a
proximal algorithm. It adds an intermediate goal-setting and search process
of exploration-exploitation and moving. Global optimization becomes a limit
case in the absence of inertia costs.

\subsection{Comparisons between local proximal and optimization algorithms}

Local search proximal algorithms were devised to
describe real life human behaviors. They work on complex 
state space and deal with radical uncertainty aspects and the physiological,
psychological, and cognitive limitations of  agents. They are also
 optimization tools.

\begin{enumerate}
\item The  question is: do these algorithms provide a realistic
description of the dynamical and stationary aspects of real life human
behaviors? The word algorithm is misleading because
these algorithms do not claim to solve optimization problems!
On the opposite, because of inertia and frictions which generate costs to
change during transitions, they help understand how humans can ultimately
reach rest points (permanent routines) which correspond to inefficient
outcomes, located far away from the optimum. Our algorithm is better viewed as a
discrete-time dynamical system describing human decision processes when
inertia matters and changing entails a cost. These algorithms involve the three
basic blocks: exploration, transition with inertia, and goal-setting blocks
with the corresponding control parameters. They allow us to describe a large
spectrum of behaviors, from muddling through behaviors, satisficing,
satisficing and \textquotedblleft worthwhile-to-move,"  to global
optimization. Attouch and Soubeyran(forthcoming)  give
several applications of  \textquotedblleft worthwhile-to-move"
behaviors and  \textquotedblleft local search proximal" algorithms. We
show how the agent can overcome inertia by using adaptive behaviors
involving long term goals and short term intermediate goals. 
Realistic decision-making models must take care of the well being
of the agents during transitions: because of inertia, agents reject
transitions with too many intermediate sacrifices (costs to change and costs
to learn how to do a new action). During transitions humans ought to survive!
Most optimization algorithms do not take  this aspect into account.

\item At each step, when using local search proximal algorithms, one has
to  solve an optimization problem. Replacing a single
optimization problem by a sequence of optimization problems
 has  advantages. In usual proximal algorithms the
quadratic cost to change  is  interpreted as a regularization
term. The first order optimality condition for
\begin{equation}  \label{Prox10}
x_{n+1} \in \mbox{argmax} \left\{ \phi(y)-\frac{1}{2\lambda} \|x_{n}-y\|^2: %
\mbox{ }y\in H \right\}.
\end{equation}
gives
\begin{equation}  \label{Prox11}
\frac{1}{\lambda} \left( x_{n+1} - x_{n} \right) - \partial \phi(x_{n+1})
\ni 0.
\end{equation}
This Eq.  (\ref{Prox11}) is an implicit discretization of the
continuous first order gradient system
\begin{equation}  \label{Prox12}
\frac{dx}{dt}(t) - \partial \phi(x(t)) \ni 0
\end{equation}
which is  the steepest ascent method. This continous dynamical system
plays a central role in optimization, differential geometry, and physics.

Proximal algorithms share most of the large time convergence properties of
this dynamical system. Convergence properties were established in
the case of a concave upper semicontinous function $\phi$  (Rockafellar, 
1976), and  in a non convex setting by using the
Kurdyka-Lojasiewicz inequality, which is valid for a large class of possibly
non-smooth functions including real analytic or semialgebraic functions
(Attouch and Bolte, 2009). These convergence
properties hold even if the initial optimization problem is ill posed with a
continuum of solutions. The algorithm asymptotically selects a particular
one,  depending on the initialization.

The definition of proximal dynamics still makes
sense in spaces without differentiable structures, where
 the norm is  replaced by some metric or relative entropy. It allows us to define
steepest ascent dynamics in metric spaces (Ambrosio et al.,
2005, for gradient flows in the space of probability measures and
applications to the Monge-Kantorovich optimal transport problem).

These algorithms also allow the general decomposition or splitting
results for structured optimization problems. This is a key property in
order to obtain convergence of best reply dynamics for potential games 
(Attouch et al. 2007, 2008).

\item Our algorithm performs better or worse than some other 
algorithms depending on the context. Let us briefly compare it  with the simulated annealing algorithm,
which is a widely used local search optimizing algorithm (a computer
context). An extensive comparison is available in Attouch and Soubeyran
(2008, working paper). Simulated annealing aims at reaching
 an optimum of an unknown function on a known and finite
state space, after a finite number of steps. Before any local
exploration, it uses a given neighbour structure of search, a given probabilistic
generation rule for new actions and a probabilistic acceptance rule, either for an
improving or, \textquotedblleft from time to time,"  a worsering action. In a
human context, local search proximal algorithms offer several improvements with realistic
features:

\begin{enumerate}
\item the state space can be infinite, and even non compact;

\item the agent does not know the geometry of the state space ex ante, hence
the content of any exploration set before doing exploration. He cannot ex ante
determine  the probability to pick a neighbor action, which forbids him
to use a probabilistic process;

\item the agent can have a semicontinuous non differentiable long term
objective;

\item the agent has not only a long term objective, but also a short term
one, which avoid him  too many temporary sacrifices during the
transition (because of costs to change): during a transition the agent ought to
survive.

\item the agent can have a more goal-oriented objective than the willingness to improve 
his situation (a
\textquotedblleft muddling through" behavior). He can set intermediate
objectives, such as  temporary satisficing (\textquotedblleft improving
enough"). This allows  us to model the  "intermediate goal-setting
processes" of Vroom (1964), and the   \textquotedblleft hard goal effect" of \ Locke and
Latham (1990).

\item the context  includes inertia through costs
to change,  attention, exploration, learning,
switching, and adaptation costs.

\item the context leads to consider learning costs 
as costs to know how a do a new action (a major case of inertia), through coercition or
good willingness.

\item the exploration process is not necessarily local. The exploration  can proceed by
visiting neighbors of neighbors.
Indeed,  an agent can start with large
steps, and end with small steps.

\item the exploration process  adapts the temporary satisficing
process to the amount of exploration which is to be done at each step. 

\item the decision-making process is rather realistic in a human
context.  The agent is not supposed to compute
sophisticated mathematical tools to take a decision.

\item the algorithm escapes from local extrema. The agent
can do large and decreasing steps at the very beginning and afterwards
from time to time. 
\end{enumerate}
\end{enumerate}

\section{Conclusion}

Our dynamical heuristic model concerns  behaviors,
defined as a successions of decisions and actions. An
action is a move in a decision space $X$. We examine why agents do
something or not, and how they choose to do it. They can
choose to do nothing, to keep their way of doing, or to innovate. The  reference
 is \textquotedblleft what and how\textquotedblright\ things where done
previously. Our model involves six interrelated blocks:

1. \textsl{Incrementalism}. The \textquotedblleft worthwhile-to-move'' or
moving block
\begin{equation}
y\in W(x)=\left\{ y\in X: \mbox{ } A(x,y)\geq \xi (x)C(x,y)\right\}
\end{equation}
\noindent with frictions and inertia models the transition from a temporary
routine $x $ to a new one $y$. The \textquotedblleft worthwhile-to-move''
principle says that a move  is  acceptable if his
estimated advantages to move $A(x,y)$ are higher than some proportion $1\geq
\xi (x)\geq 0$ of his expected costs to move. The dynamic of change follows
an acceptable transition path, where short-term sacrifices are sufficiently few.

2. \textsl{Instrumentalism}. The motivation building and goal-setting block
\begin{equation}
y\in I\left[ x,\varepsilon (x)\right] =\left\{ y\in X: \mbox{ }
g(y)-g(x)\geq \varepsilon (x)\geq 0\right\} \subset X
\end{equation}
\noindent includes improving, \textquotedblleft improving
enough,\textquotedblright   and intermediate satisficing processes with not
too much sacrificing ($g$ is an instantaneous utility function).

3. \textsl{Local exploration and search}. A basic ingredient is the local
exploration block
\begin{equation}
y\in E\left( x,r(x)\right)\subset X,
\end{equation}
\noindent with $r(x)$ being equal to the radius of the exploration set
around $x\subset X.$

4. \textsl{Heuristics}. By using parsimonious heuristics
the agent localizes, encloses the process, and cuts the cost regression
paradox (to \textquotedblleft know how to know how...''). In the context of
radical uncertainty, we do not use probabilities, but rather set membership,
which is identified by inequalities to better model the degree of
flexibility, fuzziness, and adaptability of a behavior.

5. \textsl{Punctuated dynamic}. Transitions matter, because they are
necessary to reach the  goal (a succession of jumps from a temporary
routine to another one). \textquotedblleft Walking on the
road''\ gradually becomes as important as reaching the goal itself. This
irregular dynamic articulates decisions and actions along static phases of
exploration-exploitation and dynamic moving phases.

6. \textsl{Physical, physiological, psychological, cognitive, and social
features.}  We model
many behaviors when there is a lack of knowledge (knowledge acquisition
costs) and frictions (goal-setting costs and costs of moving), (Huitt, 1992).

\smallskip

Our model links together complementary \textquotedblleft motivation
building'' (goal setting), \textquotedblleft exploration around,''  and
\textquotedblleft moving'' (changing, learning) tasks, which are the basic
blocks of the \textquotedblleft adaptive decision-making'' behaviors.
Incremental processes of decision and making are described as a succession
of comparisons between state-dependent intermediate behavioral advantages
and costs to change, with respect to intermediate advantages and costs to
stay.

We characterized  incremental behavior as a
\textquotedblleft worthwhile-to-move\textquotedblright\ behavior. At each
step, intermediate advantages to move are required to be higher than some
fraction of the costs to move to bound intermediate sacrifices. These
apparently inefficient behaviors can be rationalized as sparing behaviors.
Incremental decision-making processes (\textquotedblleft muddling
through") represent the simplest case where the motivation-and-goal process
is to do  better than before, which limits  exploration and
moving. Satisficing processes represent more goal-oriented
\textquotedblleft worthwhile-to-move" behaviors.

We rationalized incremental behaviors in a general context
including lack of knowledge (costs of exploration), lack of motivation
(costs to build or sustain motivation, goal-setting costs), inertia and
friction during the transition (costs to change, costs to learn).

We need only a metric space, there is no need of vectors.
Utility is upper-bounded to reflect limited
resources. The agent  encloses his \textquotedblleft worthwhile-to-move
relationship\textquotedblright\ in putting bounds to his control
variables, in choosing not too long an exploitation period at each step, in
using a minimum effort to move (for example when the speed of moving is not
too low and the per unit of time effort to move increases with the speed of
moving), in taking not too low a sacrifice index and putting not too heavy
weights over temporary satisfaction and disapointment (which amounts to 
\textit{high local costs to move}).

Under such assumptions,  by following a \textquotedblleft
worthwhile-to-move" process, the agent avoids intermediate sacrifices
during the transitions. This leads him to local actions in an endogenous way:
local action is no longer an hypothesis as it is for hill climbing
algorithms, it is now a consequence of inertia and frictions. Specifically,

\begin{itemize}
\item a) For a low goal-oriented behavior, where the agent just wants to
improve step by step by following the \textquotedblleft worthwhile-to-move''
relationship, the dynamical process has the local action property. It is
nested and, when the state space is complete, it converges to some final
state. Furthermore, if the agent spends a finite time for exploitation and
the speed of moving is high enough, it converges in a finite number of steps.

\item b) When the behavior is more goal-oriented,  starting from a
given initial state, the process  converges to
 a rest point, where the agent prefers to stay than to move on.
Habits and routines represent specific examples of rest points where an
agent no longer needs to think before acting.
\end{itemize}

In this context of friction, our model helps give a qualitative answer to
the question: \textquotedblleft how far from optimizing do agents behave in
real life, depending on the context\textquotedblright? Our model helps
calibrate the size of the inefficiencies (the departure of a given behavior
with respect to its substantive formulation as an optimization model). It
gives us a tool to know when the \textquotedblleft as if
hypothesis\textquotedblright\ is indicated, when optimization is a good
enough approximation of a given behavior. To this purpose our model must

i) specify how to catch behavioral inefficiencies, how to model lack of
knowledge, frictions and goal-setting inefficiencies and efficiencies. This
problem becomes more complicated when the model involves costs to move;

ii) define dynamic inefficiency indices, which amounts to calibrate the
inefficiency gaps of a behavior with respect to its substantive formulation;

iii) link the inefficiency gaps of a behavior to the characteristics of both
the agent and his environment.

Concerning costs to move, inefficiencies indices also include the total
costs to move, the mean velocity, and the intermediate sacrifices. Several
quality-cost ratios were defined.

At a mathematical level, this led us to revisit the \textquotedblleft as if
hypothesis'' and put to the fore the \textquotedblleft local search and
proximal algorithms'' which mix local search and proximal algorithms.

\section{References}

Abel, A. (1990). Asset prices under habit formation and catching up with the
Joneses. \textit{American Economic Review}, 80 (2), 38-42.

\vskip 1mm

\noindent Ambrosio, L., Gigli, N., and Savare, G. (2005). Gradient flows in
metric spaces and in the space of probability measures. \textit{Lectures in
Mathematics ETH Z$\ddot{u}$rich}, Birkh$\ddot{a}$user Verlag.

\vskip 1mm

\noindent Aarts, E., and Lenstra, K. (2003). Local search. \textit{Princeton
University Press}.

\vskip 1mm

\noindent Attouch, H., and Bolte, J. (2009). On the convergence of the
proximal algorithm for non smooth functions involving analytic features. 
\textit{Math. Program.}, Ser. B, 116, (1-2), 5-16.

\vskip 1mm

\noindent Attouch, H., Bolte, J., Redont, P., and Soubeyran, A. (2008).
Alternating proximal algorithms for weakly coupled convex minimization
problems. Applications to dynamical games and PDE's. \textit{J. 
Convex Anal.}, 15, (3), 485-506.

\vskip 1mm

\noindent Attouch, H., Goudou, X., and Redont P. (2000). The heavy ball with
friction method. \textit{Communication in Contemporary Mathematics}, 2,
(1), 1-34.

\vskip1mm

\noindent Attouch, H., Redont, P., and Soubeyran, A. (2007). A new class of
alternating proximal minimization algorithms with costs-to-move. 
\textit{SIAM J.  Optim.}, 18, (3), 1061-1081.

\vskip 1mm

\noindent Attouch, H., Soubeyran A. (2006). Inertia and reactivity in
decision making as cognitive variational inequalities'', \textit{J. 
Convex Anal.}, 13, (2), 207-224.

\vskip 1mm

\noindent Attouch H., and Soubeyran, A. (2006). A cognitive approach of the
Ekeland theorem, section 3.4. of \textquotedblleft Variational analysis in
Sobolev and BV spaces: Applications to PDE and Optimization'', Editors H.
Attouch - G. Buttazzo - G. Michaille, \textit{MPS-SIAM series on Optimization%
}, MP06.

\vskip 1mm

\noindent Attouch, H., Soubeyran A. (2008). Comparison between
 local proximal and  simulated
annealing algorithms, Working paper, GREQAM.

\vskip 1mm

\noindent Attouch, H., and Teboulle M. (2004).  A
Regularized Lotka-Volterra dynamical system as a continuous proximal-like
method in optimization, \textit{J. Optim. Theory Appl}., 121, (3),
541-580.

\vskip 1mm

\noindent Aubin, J.P., and Ekeland, I. (1983). \textit{Applied non
linear analysis}, Wiley Interscience.

\vskip 1mm

\noindent Aubin, J.P. (2005).  Evolution tychastique,
stochastique et contingente, mimeo.

\vskip 1mm

\noindent Bettman, J., Luce M., and Payne, J. (1998). Constructive consumer
choice process. \textit{Journal of Consumer Research}, 25, 187-217.

\vskip 1mm

\noindent Busemeyer, J.R., and Diederich, A. (2002). Survey of decision field
theory. \textit{Math. Soc. Sci.}, 43, (3), 345-370.

\vskip 1mm

\noindent Carroll, C. (2001). Solving consumption models with multiplicative
habits. \textit{Economics Letters}, 68, (1), 67-77.

\vskip 1mm

\noindent Conlisk, J. (1996). Why bounded rationality? \textit{Journal of
Economic Literature}, 34, 669-700.

\vskip 1mm

\noindent Ekeland, I. (1974). On the variational principle. \textit{J. Math.
Anal. Appl.}, 47, 325-353.

\vskip 1mm

\noindent Fan K. (1972). A minimax inequality and applications, in \textit{%
Inequalities} 3, O Sisha Ed, 103, 113, Academic Press.

\vskip 1mm

\noindent Franklin, A. (1978). Consumer decision making and health
maintenance organizations: A review. \textit{Medical Care}, 16 , 1-13.

\vskip 1mm

\noindent Friedman J. (1953). The methodology of positive economics, \textit{%
Essays in positive economics}, University Press of Chicago.

\vskip 1mm

\noindent Gigerenzer, G., and Todd, P. (1999). Simple heuristics that make
us smart. Oxford University Press.

\vskip 1mm

\noindent Gilovich, T., Griffin, D., and Kahneman D. (2002). Heuristics and
biases, Cambridge University Press.

\vskip 1mm

\noindent Goudou, X., and Munier, J. (2005). Asymptotic behavior of
solutions of a gradient like integrodifferential Volterra inclusion. \textit{%
Advances in Math. Sciences and Applic.}, 15, 509-525.

\vskip 1mm

\noindent Huitt, H. (1992). Problem solving and decision making:
Consideration of individual differences using the Myers-Briggs type
indicator, \textit{Journal of Pyschological Type}, 24, 33-44.

\vskip 1mm

\noindent Iusem, A. (1995). Proximal point methods in optimization.
Instituto de Matematica Pura e Applicada, Rio de Janeiro, Brasil.

\vskip 1mm

\noindent Klemperer, P. (1995). Competition when consumers have switching
costs: An overview with applications to industrial organization,
macroeconomics, and international trade. \textit{Review of Economic Studies}%
, 62, 515-539.

\vskip 1mm

\noindent Lindblom, C. (1959). The science of muddling through. \textit{%
Public Administrative Review}, 19, 79-99.

\vskip 1mm

\noindent Lipman, B., and Wang, R. (2000). Switching costs in frequently
repeated games. \textit{Journal of Economic Theory}, 93, 149-190.

\vskip 1mm

\noindent Lipman, B., and Wang, R. (2006). Switching costs in infinitely
repeated gmes. Queen's Economics Department, Boston University, Working
paper 1032.

\vskip 1mm

\noindent Locke, E., and Latham, G. (1990). A theory of goal setting and
task performance, Englewoods Cliffs, NJ: Prentice Hall.

\vskip 1mm

\noindent Martinez-Legaz, JE., and Soubeyran, A. (2002). Learning from
errors. UAB, Working paper LEA, University of Barcelona.

\vskip 1mm

\noindent Martinez-Legaz, JE., and Soubeyran, A. (2007). A Tabu search
scheme for astract problems, with applications to the computation of fixed
points. \textit{J. Math. Anal. Appl.}, 73, (1), 1-8.

\vskip 1mm

\noindent Nelson, R., and Winter, S. (1997). An evolutionary theory of
economic change, in \textit{Resources Firms and Strategies}, Oxford
Management Readers, Oxford.

\vskip 1mm

\noindent O' Donoghue, T., and Rabin, M. (1999). Doing it now or later. 
\textit{\ American Economic Review}, 89, (1), 103-124.

\vskip 1mm

\noindent Payne, J., Bettman, J., and Johnson E. (1988). 
Adaptive strategy selection in decision making, \textit{Journal of
Experimental Psychology, Learning, Memory and Cognition}, 14, (3), 534-553.

\vskip 1mm

\noindent Rockafellar, R.T. (1976). Monotone operators and the proximal
point algorithm. \textit{SIAM J. Control Opt.}, 14, (5), 877-898.

\vskip 1mm

\noindent Rockafellar, R.T., and Wets, R. (2004). 
\textit{Variational Analysis}, Springer- Verlag.

\vskip 1mm

\noindent Rumelt, R. (1990). Inertia and Transformation. mimeo.

\vskip 1mm

\noindent Selten, R. (1998). Features of experimentally observed bounded
rationality. \textit{European Economic Review}, 42, 413-436.

\vskip 1mm

\noindent Selten, R. (1998. Aspiration adaptation theory. \textit{J.
Math. Psychology}, 42, 191-214.

\vskip 1mm

\noindent Simon, H. (1955). A behavioral model of rational choice. \textit{%
Quaterly Journal of Economics}, 69, 99-118.

\vskip 1mm

\noindent Simon, H. (1967). From substantive to procedural rationality, in 
\textit{Methods and Appraisal in Economics}, Latsis editor.

\vskip 1mm

\noindent Simon, H. (1982). The science of the artificial. 2nd Edition,
Cambridge, MA: MIT Press.

\vskip 1mm

\noindent Simon, H. (1987. Satisficing, in \textit{the New Palgrave: A
dictionary of Economics}, Eds Eatwell J., Milgate M., Newman P., London, Mac
Millan, 243-245.

\vskip 1mm

\noindent Sinclair-Desgagn\'{e}, B., and Soubeyran, A. (2000). A theory of
routines as mindsavers. Cirano, Working Paper, 2000, s52.

\vskip 1mm

\noindent Singh, S., Watson, B.,and Srivastava, P. (1997). \textit{Fixed
point theory and best approximation: The KKM map principle}. Kluwer.

\vskip 1mm

\noindent Sobel, J. (2000). Economists' models of learning. \textit{J.
Econom. Theory}, 94 , 241-261.

\vskip 1mm

\noindent Soubeyran, A. (2006). Adaptative satisficing processes: How to
take benefit of both negative and positive knowledge acquisition.

\vskip 1mm

\noindent Soubeyran, A., and Soubeyran, B. (2009a). Variational rationality as
changing organizational capabilities: Worthwhile transitions between
temporary routines. Working paper GREQAM.

\vskip 1mm

\noindent  Soubeyran, A., and Soubeyran, B. (2009b). A value creation-value
appropriation theory of organizational change: A variational
approach. Working paper GREQAM.

\vskip 1mm

\noindent Tyson, C. (2003).  Revealed preference analysis of boundedly rational choice. PhD
thesis, Stanford University.

\vskip 1mm

\noindent Vroom, VH. (1964). Work and motivation. New York, Wiley.

\vskip 1mm

\noindent Williamson, O. (1975). Market and hierarchies: Analysis and anti
trust implications. NY: Free Press.

\end{document}